\documentclass[12pt,leqno]{article}
\usepackage{oxford3}
\usepackage{graphicx}
\usepackage{amsfonts}
\usepackage{amsmath}
\usepackage{rotating}
\usepackage{enumerate}
\usepackage{amssymb, latexsym,amsthm}
\usepackage{multirow,array,epsfig}
\usepackage{verbatim}
\usepackage{color}
\usepackage{mathtools}
\usepackage{booktabs}
\usepackage{caption}
\usepackage{url}
\usepackage{accents}

\theoremstyle{definition}

\newtheorem{assumption}{Assumption}[section]
\newtheorem{definition}{Definition}[section]
\newtheorem{remark}{Remark}[section]

\newtheorem{procedure}{Procedure}[section]

\theoremstyle{theorem}
\newtheorem{theorem}{Theorem}[section]

\newtheorem{corollary}{Corollary}[section]

\numberwithin{equation}{section}
\renewcommand{\baselinestretch}{1.05}
\newcommand\bGa{\mbox{\boldmath${\Gamma}$}}

\newcommand\bXi{\mbox{\boldmath${\Xi}$}}

\newcommand\bff{{\bf f}}

\newcommand\mbH{{\mathbb H}}

\newcommand\bM{{\bf M}}

\newcommand\bQ{{\bf Q}}
\newcommand\bR{{\bf R}}

\newcommand\mbR{{\mathbb R}}
\newcommand\bs{{\bf s}}

\newcommand\bW{{\bf W}}

\newcommand\bX{\boldsymbol{X}}

\newcommand\bZ{{\bf Z}}

\newcommand\mcN{{\mathcal N}}

\newcommand\cS{{\mathcal S}}

\newcommand\cU{{\mathcal U}}

\newcommand\tp{{t^\prime}}

\DeclareMathOperator{\E}{E}

\DeclareMathOperator{\Span}{span}
\DeclareMathOperator{\tr}{tr}
\DeclareMathOperator{\Tr}{Tr}
\DeclareMathOperator{\Var}{Var}

\newcommand{\notate}{\textcolor{black}}

\begin{document}

\title{Testing Separability of Functional Time Series}
\author{Panayiotis Constantinou\\
{\small Pennsylvania State University}
\and
Piotr Kokoszka\\
{\small Colorado State University}
\and
Matthew Reimherr\thanks{\small {\em Corresponding author:}
Department of Statistics,
Pennsylvania State University, 411 Thomas Building,
University Park, PA 16802, USA. \
mreimherr@psu.edu \
(814) 865-2544 } \\
{\small Pennsylvania State University}
}
\date{}
\maketitle

\begin{abstract}
We derive and study a significance test for determining if a panel of
functional time series is separable.  In the context of this paper,
separability means that the covariance structure factors into the
product of two functions, one depending only on time and the other
depending only on the coordinates of the panel.  Separability is a
property which can dramatically improve computational efficiency by
substantially reducing model complexity.  It is especially useful for
functional data as it implies that the functional principal components
are the same for each member of the panel.  However such an assumption
must be verified before proceeding with further inference. Our
approach is based on functional norm differences and provides a test
with well controlled size and high power.
\notate{We establish our procedure quite generally, allowing one to test separability of autocovariances as well.} In addition to an
asymptotic justification, our methodology is validated by a simulation
study. It is applied to functional panels of particulate pollution and
stock market data.
\end{abstract}

\section{Introduction} \label{s:i}
Suppose $\left\{ X(\bs, t) , \bs \in \mathbb{R}^2, t \in \mathbb{R}
\right\}$ is a real--valued spatio--temporal random field,
with the coordinate $\bs$ referring to space,  and $t$ to time.
The field $X(\cdot, \cdot)$ is said to be {\em separable} if
\[
\text{Cov} (X(\bs_1, t_1), X(\bs_2, t_2))
= u(\bs_1, \bs_2) v(t_1, t_2),
\]
where $u$ and $v$ are, respectively, spatial and temporal covariance
functions.  Separability is discussed in many textbooks, e.g
\citetext{cressie:wikle:2015}, Chapter 6. It has been extensively used
in spatio--temporal statistics because it leads to theoretically
tractable models and computationally feasible procedures; some recent
references are \citetext{hoff:2011}, \citetext{paul:peng:2011},
\citetext{sun:li:genton:2012}.  Before separability
is assumed for the reasons noted above, it must be
tested. Tests of separability are reviewed in
Mitchell {\it et al.}  (\citeyear{mitchell:genton:gumpertz:2005},
\citeyear{mitchell:genton:gumpertz:2006}) and \citetext{fuentes:2006}.

Time series of weather or pollution
related  measurements obtained at spatial locations typically
exhibit strong periodic patterns. An approach to accommodate
this periodicity  is
to divide the time series of such type  into segments, each segment
corresponding to a natural period. For example, a  periodic time series
of maximum daily temperatures at some location can be viewed
as a stationary time series of functions, with one function per year. If the
measurements are available at many locations $\bs_k$, this
gives rise to a data structure of the form
\[
X_n(\bs_k; t_i), \ k =1, \ldots, S, \ i= 1, \ldots, I (=365), \ n = 1, \ldots, N,
\]
where $n$ indexes year, and $t_i$ the day within a year. Time series
of functions are discussed in several books, e.g. \citetext{bosq:2000},
\citetext{HKbook}, \citetext{KRbook}, but research on spatial fields
or panels of time series of functions is relatively new, e.g.
\citetext{kokoszka:reimherr:woelfing:2016},
Gromenko {\it et al.}
(\citeyear{gromenko:kokoszka:reimherr:2016},
\citeyear{gromenko:kokoszka:sojka:2017}),
\citetext{french:KSH:2016},
\citetext{tupper:matteson:anderson:2017},
\citetext{liu:ray:hooker:2017} and
\citetext{shang:hyndman:2017}.
Testing
separability of spatio--temporal functional data of the above
form is investigated in \citetext{constantinou:2017}, 
\citetext{aston2017tests} and \citetext{bagchi:dette:2017},  under the
assumption that the fields $X_n(\cdot, \cdot), 1 \le n \le N,$ are
independent. No tests are currently
available for testing separability in the presence
of temporal dependence across $n$.
In a broader setting, the data that motivate this research have the
form of functional panels:
\begin{equation} \label{e:panel}
\boldsymbol{X}_n(t) = [X_{n1}(t), X_{n2}(t), ... , X_{nS}(t)]^T,
\quad     1 \le n \le N. \\
\end{equation}
Each $X_{ns}(\cdot )$ is a curve, and all curves are defined on the
same time interval. The index $n$ typically stands for day, week,
month or year.  For instance, $X_{ns}(t)$, can be the exchange rate
(against the Euro or the US Dollar) of currency $s$ at minute $t$ of
the $n$th trading day, or $X_{ns}(t)$ can be the stock price of
company $s$ at minute $t$ of the $n$th trading day. Another
extensively studied example is daily or monthly yield curves for a
panel of countries, e.g.  \citetext{ang:bekaert:2002},
\citetext{bowsher:meeks:2008}, \citetext{hays:shen:huang:2012},
\citetext{kowal:matteson:ruppert:2016}, among others.  As for scalar
data, the assumption of separability has numerous benefits including a
simpler covariance structure, increased estimation accuracy, and
faster computational times. In addition, in the contexts of functional
time series, separability implies that the optimal functions used for
temporal dimension reduction are the same for each member (coordinate)
of the panel; information can then be pooled across the coordinates to
get better estimates of these functions. We elaborate on this point in
the following. However, if separability is incorrectly assumed, it
leads to serious biases and misleading conclusions. A significance
test, which accounts for the temporal dependence present in all
examples listed above, is therefore called for.  The derivation of
such a test, and the examination of its properties, is the purpose of
this work.  \notate{Our procedure is also applicable to testing
separability of the autocovariance at any lag. }  We
will demonstrate that it works well in situations where the tests of
\citetext{constantinou:2017} and \citetext{aston2017tests} fail.

The remainder of the paper is organized as follows. In Section
\ref{s:ap},  we formulate the assumptions, the definitions, and the
problem. In Section \ref{s:da}, we derive the test and provide the
required asymptotic theory.  Section \ref{s:det} focuses on
details of the implementation. In Section~\ref{s:sim},  we present
results of a simulation study, and, finally, in Section~\ref{s:appl} we
apply our procedure to functional panels of Nitrogen Dioxide levels on
the east coast of the United States and to US stock market data.

\section{Assumptions and problem formulation} \label{s:ap} We assume
that the $\bX_n$ in \eqref{e:panel} form a strictly stationary
functional time series of dimension $S$. To simplify notation, we
assume that all functions are defined on the unit interval $[0,1]$
(integrals without limits indicate integration over $[0,1]$).  We
assume that they are square integrable in the sense that $ E \| X_{ns}
\|^2 = E \int X_{ns}^2(t) dt < \infty.  $ Stationarity implies that
the \textit{lagged} covariance function can be expressed as
\[
\text{Cov}(X_{ns}(t),X_{{n+h},s'}(t'))=c^{(h)}(s,t,s',t').
\]
\notate{
We aim to test the null hypothesis
\begin{equation} \label{e:H_0}
H_0:\ c^{(h)}(s,t,s',t')=c_1^{(h)}(s,s') c_2^{(h)}(t,t'),
\quad  s,s' \in \{1,2,\ldots, S\}; \ t,t' \in [0,1],
\end{equation}
for a fixed value of $h$. The most important setting is when $h = 0$,
i.e., testing separability of the covariance function, but other lags
can be considered as well.  }

To derive the asymptotic distribution of our  test statistic we impose
a weak dependence condition on the $\bX_n$. We use the
concept of $L^p$--$m$--approximability introduced in
\citetext{hormann:kokoszka:2010}, see also Chapter 16
of \citetext{HKbook}.  Suppose $\mbH$ is a separable Hilbert
space. Let
$p\geq 1$ and let $L_\mbH^p$ be the space of $\mbH$--valued random
elements $X$ such that
\[
\nu_p(X) = \big( E\|X\|^p \big)^{1/p} <\infty.
\]

\begin{definition}\label{d:approx}
The sequence, $\{Z_n\}, \ Z_n \in L_\mbH^p $,
is {\em  $L^p$--$m$--approximable} if the following
conditions hold:
\begin{enumerate} \itemsep0em
\item There exists a sequence
$\left\{ u_n \right\}$ of iid elements in an abstract measurable space $\cU$
such that
\[
Z_n = f(u_n, u_{n-1}, \dots),
\]
for a measurable function  $f: \cU^{\infty} \to \mbH$.
\item For each integer $M> 0$,  consider an approximating
sequence $Z_{n,M}$ defined by
\[
Z_{n, M} = f(u_n, u_{n-1}, \ldots, u _{n-M},
u^\star_{n-M-1}, u^\star_{n-M-2},\dots),
\]
where the sequences $\left\{ u^\star_n \right\} = \left\{ u^\star_n(n,m) \right\}$
are copies of
$\left\{ u_n \right\}$  independent across $m$ and $n$ and independent of
the original sequence $\left\{ u_n \right\}$.
We assume that $Z_{n,M}$ well approximates $Z_n$ in the sense that
\begin{equation} \label{e:approx}
\sum_{M=1}^\infty  \nu_p (Z_n - Z_{n, M })
< \infty.
\end{equation}
\end{enumerate}
\end{definition}

Condition 1 of Definition \ref{d:approx} implies that the
sequence is strictly stationarity and ergodic. The essence of
Condition 2 is that the dependence of $f$ on the innovations far in
the past decays so fast that these innovations can be replaced by
their independent copies.  Such a replacement is asymptotically
negligible in the sense quantified by \eqref{e:approx}. Similar
conditions, which replace the more restrictive
assumption of a linear moving average
with summability conditions on its
coefficients,  have been used for at least a decade, see e.g.
\citetext{shao:wu:2007} and references therein. We work with
Definition~\ref{d:approx} as it is satisfied by most time series
models, including functional time series,
 and provides a number of desirable asymptotic properties
including the central limit theorem,
\textcolor{black}{see Chapter 16 of
\citetext{HKbook} and \citetext{kokoszka:reimherr:2013}, among many
 other references.  The conditions in Definition~\ref{d:approx} cannot
be verified, they are analogous to mixing or summability of cumulants
conditions which have been imposed in theoretical time series analysis
research.} We therefore make the following
assumption.

\begin{assumption}\label{a:main}
The $\bX_n$ form an  $L^4$--$m$--approximable sequence
 in $\mbH= (L^2([0,1]))^S$.
\end{assumption}

We use
tensor notation analogous to \citetext{aston2017tests}.
Let $\mbH_1$ and $\mbH_2$ denote two real separable Hilbert spaces
with bases $\{u_i\}$ and $\{v_j\}$,  respectively. We define
$\mbH=\mbH_1\otimes\mbH_2$ to be the tensor product Hilbert space.
The tensors $\{u_i\otimes v_j\}$ form a basis for $\mbH$. In
other words, the tensor product Hilbert space can be obtained by
completing of the set $\Span\{ u_i\otimes v_j: i=1,\dots \ j = 1,
\dots\}$, under the following inner product:
\[
\langle u_i\otimes v_j, u_k\otimes v_\ell \rangle = \langle u_i, u_k\rangle \langle v_j, v_\ell\rangle,\quad u_i, u_k \in \mbH_1, v_j, v_\ell \in \mbH_2.
\]
In the context of our study $\mbH_1=\mbR^{S}$ and
$\mbH_2=L^2([0,1])$.
Therefore the tensor product Hilbert space in our context is
$\mbH=\mbH_1\otimes\mbH_2=\mbR^{S}\otimes
L^2([0,1])=(L^2([0,1]))^S=:L_2^S$, where we omit $[0,1]$ for
simplicity.  Each $\boldsymbol{X}_n$ is thus an element of a tensor
space, formed by the tensor product between two real separable Hilbert
spaces, $\boldsymbol{X}_n \in \mbH_1 \otimes \mbH_2$. We denote by
$\cS(\mbH_1 \otimes \mbH_2)$ the space of Hilbert-Schmidt operators
acting on $\mbH_1 \otimes \mbH_2$. Note that $\{u_i\otimes v_j\otimes
u_k\otimes v_\ell\}$ is a basis for $\cS(\mbH_1 \otimes \mbH_2)$.  The
covariance operator \textcolor{black}{between $\bX_n$ and  $\bX_{n+h}
  \in\mbH=\mbH_1 \otimes \mbH_2$}, $C^{(h)} = \E [\bX_n \otimes
\bX_{n+h}] \in \cS(\mbH_1 \otimes \mbH_2)$, is called separable if
\begin{equation} \label{e:H0-C}
 C^{(h)} = C_1^{(h)} \widetilde \otimes C_2^{(h)},
\end{equation}
where $C_1^{(h)}$ is a covariance operator over $\mbH_1$ and $C_2^{(h)}$ is a
covariance operator over $\mbH_2$. We define $C_1^{(h)} \widetilde \otimes
C_2^{(h)}$ as a linear operator on $\mbH_1 \otimes \mbH_2$ satisfying
\[
(C_1^{(h)} \widetilde \otimes C_2^{(h)})( u \otimes v) = (C_1^{(h)} u) \otimes (C_2^{(h)} v),
\quad \forall u \in \mbH_1,\forall v \in \mbH_2.
\]
The covariance operator \textcolor{black}{between $\bX_n$ and  $\bX_{n+h}
\in L_2^S$ is in  $\cS(L_2^S)$, i.e. it is an
integral operator with the kernel $c^{(h)}$.}  Relation \eqref{e:H0-C} is
then equivalent to $H_0$ stated as \eqref{e:H_0}  above.

\section{Derivation of the test and its asymptotic justification}
\label{s:da} To test hypothesis \eqref{e:H0-C}, we propose a statistic
which quantifies the difference between $\widehat C_1^{(h)} \widetilde
\otimes \widehat C_2^{(h)} $ and $\widehat C^{(h)}$:
\begin{equation} \label{e:test}
\widehat T
= N \| \widehat C_1^{(h)} \widetilde\otimes \widehat C_2^{(h)}
- \widehat C^{(h)} \|_{\cS}^2,
\end{equation}
where $\widehat C_1^{(h)}, \widehat C_2^{(h)}, \widehat C^{(h)}$
are estimates defined
below, and $\|\cdot\|_{\cS}$ is the Hilbert-Schmidt norm.
\textcolor{black}{
The statistic \eqref{e:test} is
a normalized distance between the estimator valid
under the restriction imposed by $H_0$ and a general unrestricted
estimator. The term
$\widehat C_1^{(h)} \widetilde\otimes \widehat C_2^{(h)}$
is an estimator of the product
$c_1^{(h)}(\cdot, \cdot) c_2^{(h)}(\cdot, \cdot)$ in
 \eqref{e:H_0} (the autocovariance under separability),
whereas  $\widehat C^{(h)}$ is an estimator of
the unrestricted spatio--temporal autocovariance function
$c^{(h)}(\cdot, \cdot, \cdot, \cdot)$. While  $\widehat C^{(h)}$
is not difficult to define, it is not
obvious how to define $\widehat C_1^{(h)}$ and $\widehat C_2^{(h)}$.
This section explains how we define the estimators in \eqref{e:test}
and what their joint asymptotic distribution is. This will allow us
to derive the asymptotic properties of $\widehat T$.
}%textcolor

The
asymptotic null distribution involves the covariance operator of
$\widehat C_1^{(h)} \widetilde\otimes \widehat C_2^{(h)} - \widehat C^{(h)}$, which we
denote by $\bQ^{(h)}$. Note that $\bQ^{(h)} \in \cS(\cS(\mbH_1 \otimes \mbH_2))$,
i.e. it is an operator acting on $\cS(\mbH_1 \otimes \mbH_2)$.  Therefore,
it can be expanded using the basis functions of the
form $\{u_i\otimes v_j\otimes u_k\otimes v_\ell\otimes u_m\otimes
v_n\otimes u_p\otimes v_q\}$. In the context of \eqref{e:panel},  $\bQ^{(h)} \in
\cS(\cS(L_2^S))$.

We now define the estimators appearing in \eqref{e:test} and obtain
their limiting behavior even in the case where $C^{(h)}$ is not separable.
A natural estimator for the general covariance, $C^{(h)}$, is given by\\
\[
\widehat{C}^{(h)}=\frac{1}{N-h} \sum_{n=1}^{N-h}
(\boldsymbol{X}_n - \boldsymbol{\hat{\mu}})
\otimes(\boldsymbol{X}_{n+h} - \boldsymbol{\hat{\mu}})
\in \cS(L_2^S),
\]
where $ \boldsymbol{X}_n(t) = [X_{n1}(t), X_{n2}(t), ... ,
X_{nS}(t)]^T$ , $1 \le n \le N,$ and
$\boldsymbol{\hat{\mu}}(t)=[\hat{\mu}_{1}(t), \hat{\mu}_{2}(t), ... ,
\hat{\mu}_{S}(t)]^T$ with
$\hat{\mu}_{s}(t)=\frac{1}{N}\sum_{n=1}^{N}X_{ns}(t)$, $1 \le s \le
S$.
Since centering by the sample mean is asymptotically negligible, we
assume, without loss of generality and to ease
the notation, that our data are centered,  so the
estimator takes the form
\begin{equation} \label{e:hat-Ch}
\widehat{C}^{(h)}=\frac{1}{N-h} \sum_{n=1}^{N-h} \boldsymbol{X}_n\otimes \boldsymbol{X}_{n+h},
\end{equation}
equivalently, the kernel of $\widehat{C}^{(h)}$ is
\[
\hat{c}^{(h)}(s,t,s',t')=\frac{1}{N-h} \sum_{n=1}^{N-h} X_{ns}(t)X_{{n+h},s'}(t').
\]
Under $H_0$, $C^{(h)}=C_1^{(h)}\widetilde\otimes C_2^{(h)}$ with $C_1^{(h)}\in\cS(\mbH_1) = \cS(
\mbR^{S})$, $C_2^{(h)}\in\cS(\mbH_2) = \cS(L^2([0,1]))$ and $C^{(h)}
\in \cS(\mbH)=\cS(\mbH_1 \otimes \mbH_2) =\cS(L_2^S).$
To obtain estimators  for $C_1^{(h)}$ and $C_2^{(h)}$, we
utilize the trace and the partial trace operators. For any trace-class
operator $T$, see e.g. Section 13.5 of \citetext{HKbook} or
Section 4.5 of \citetext{hsing:eubank:2015}, its trace is defined by
\[
\Tr(T):=\sum_{i = 1}^\infty \langle   Te_i , e_i \rangle,
\]
where $(e_i)_{i\ge 1}$ is an orthonormal basis. It is invariant with
respect to the basis.  The partial-trace operators are defined
as
\[
 \Tr_1(A \widetilde\otimes B)=\Tr(A)B, \ \ \
A \in \mbH_1,\ B  \in \mbH_2,
\]
and
\[
 \Tr_2(A \widetilde\otimes B) = \Tr(B)A, \ \ \
 A \in \mbH_1,\   B \in \mbH_2.
\]
This means that $\Tr_1$ and $\Tr_2$ are bilinear forms that satisfy $
\Tr_1 : \mbH_1\otimes \mbH_2 \rightarrow \mbH_2 $ and $\Tr_2 :
\mbH_1\otimes \mbH_2 \rightarrow \mbH_1$.  In general, the trace of
any element of $T \in \mbH_1 \otimes \mbH_2$ can be defined using
proper basis expansions.  More specifically, let $u_1,u_2, \dots$ be
an orthonormal basis for $\mathbb{H}_1$ and $v_1,v_2, \dots$ an
orthonormal basis for $\mathbb{H}_2$. Then a basis for
$\mathbb{H}_1\otimes \mathbb{H}_2$ is given by $\{u_i \otimes v_j: i =
1,2,\dots , j = 1,2,\dots\}$. Let $T:\mathbb{H}_1\otimes
\mathbb{H}_2\rightarrow\mathbb{H}_1\otimes \mathbb{H}_2$. Then,  the
trace of $T$ is defined by:
\[
\Tr(T)=\sum_{i \ge1}\sum_{j \ge1} \langle T(u_i \otimes v_j),
 u_i \otimes v_j \rangle,\quad
\Tr : \mathbb{H}_1\otimes \mathbb{H}_2 \rightarrow \mbR.
\]
If $T=A \widetilde\otimes B$,
 the partial-trace operators in terms of a basis are defined as
\begin{align*}
\Tr_1(T)&=\Tr_1(A \widetilde\otimes B)=\Tr(A)B=\sum_{i \ge 1} \langle Au_i, u_i \rangle B\\
&=\sum_{i \ge 1} \langle Au_i, u_i \rangle \sum_{j \ge 1} B_j v_j , \quad \forall A \in \mathbb{H}_1,\quad \forall  B \in \mathbb{H}_2,
\end{align*}
and
\begin{align*}
\Tr_2(T)&=\Tr_2(A \widetilde\otimes B) = \Tr(B)A=\sum_{j \ge 1} \langle Bv_j, v_j \rangle A\\
&=\sum_{j \ge 1} \langle Bv_j, v_j \rangle \sum_{i \ge 1} A_i u_i,  \quad \forall A \in \mathbb{H}_1,\quad \forall B \in \mathbb{H}_2.
\end{align*}
In the context of functional panels, let $u_1,u_2, ..., u_S$ be an
orthonormal basis for $\mbR^{S}$ and $v_1,v_2, ...$ an orthonormal
basis for $L^2([0,1])$. Then a basis for $L_2^S$ is given by $\{u_i
\otimes v_j: i = 1,2,\dots, S, j = 1,2,\dots\}$. Recall that the
products $u_i \otimes u_{k}$, viewed as operators, form a basis for
$\cS(\mbR^{S})$, that is a basis for the space of Hilbert-Schmidt
operators acting on $\mbR^{S}$.  Similarly $\{v_j \otimes v_{\ell}\}$
is a basis for $\cS(L^2([0,1]))$.  Finally $\{u_i \otimes v_j \otimes
u_{k}\otimes v_{\ell}\}$ is a basis for $\cS(L_2^S)$.  The basis
expansion of $C^{(h)}$ is given by
\[
 \sum_i \sum_j \sum_{k}\sum_{\ell}C_{ijk\ell}^{(h)}  u_i \otimes v_j \otimes u_{k}\otimes v_{\ell}.
\]
Therefore its trace is given by
\[
\Tr(C^{(h)})=\sum_i \sum_j C_{ijij}^{(h)}.
\]
Under the assumption of separability, i.e. $C^{(h)}=C_1^{(h)}\widetilde\otimes C_2^{(h)}$, the partial trace with respect to $\mathbb{H}_1$ in terms of a basis is given by
\[
\Tr_1(C^{(h)}) =\Tr_1(C_1^{(h)}\widetilde\otimes C_2^{(h)})= \Tr(C_1^{(h)}) C_2^{(h)} = \sum_j \sum_{\ell} (\sum_i C_{iji\ell}^{(h)}) v_j \otimes v_{\ell}\quad \text{with} \quad C_{2,j\ell}^{(h)} = \sum_i C_{iji\ell}^{(h)},
\]
and with respect to $\mathbb{H}_2$ is given by
\[
\Tr_2(C^{(h)}) =\Tr_1(C_1^{(h)}\widetilde\otimes C_2^{(h)}) = \Tr(C_2^{(h)}) C_1^{(h)} = \sum_i \sum_{k} (\sum_j C_{ijkj}^{(h)}) u_i \otimes u_{k}\quad \text{with}\quad C_{1,ik}^{(h)} = \sum_j C_{ijkj}^{(h)}.
\]
Under the assumption of separability,
we define estimators of $C_1^{(h)}$ and $C_2^{(h)}$ as
\begin{equation} \label{e:hat-C12}
\widehat{C}_1^{(h)} = \frac {1} {\Tr(\widehat{C}^{(h)})}\Tr_2(\widehat{C}^{(h)})
\quad  \text{and} \quad \widehat{C}_2^{(h)}
= \Tr_1(\widehat{C}^{(h)}),
\end{equation}
where $\widehat{C}_1^{(h)}$ is an $S \times S$ matrix and $\widehat{C}_2^{(h)}$ is
a temporal covariance operator. The intuition behind the above
estimators is that $\Tr(C^{(h)})C^{(h)}=\Tr_2(C^{(h)})\widetilde\otimes\Tr_1(C^{(h)})$. Note
that the decomposition $C^{(h)}=C_1^{(h)}\widetilde\otimes C_2^{(h)}$ is not unique
since $C_1^{(h)}\widetilde\otimes C_2^{(h)}=(\alpha
C_1^{(h)})\widetilde\otimes(\alpha^{-1}C_2^{(h)})$ for any $\alpha\ne0$, however
the product $C_1^{(h)}\widetilde\otimes C_2^{(h)}$ is.

To derive the asymptotic distribution of the test statistic $\widehat
T$ defined in \eqref{e:test}, we must first derive the joint
asymptotic distribution of $\widehat C^{(h)}, \widehat C_1^{(h)},
\widehat C_2^{(h)}$.  A similar strategy was used in
\citetext{constantinou:2017}.  However, there the observations were
assumed to be independent and more traditional likelihood methods were
used to derive the asymptotic distributions.  Here, we take a
different approach, instead using the CLT for $\widehat C^{(h)}$, and
then leveraging a Taylor expansion over Hilbert spaces to obtain the
joint asymptotic distribution of $\widehat C^{(h)}, \widehat
C_1^{(h)}, \widehat C_2^{(h)}$.  In this way, we are able to relax
both the independence and Gaussian assumptions from
\citetext{constantinou:2017}.  The result is provided in
Theorem~\ref{thm:1}. Due to the temporal dependence, the covariance
operator of the limit normal distribution is a suitably defined
long--run covariance operator.  It has a very complex, but explicit
and computable, form, which is displayed in  Supporting
Information, where all theorems that follow are also proven.

Recall that we are interested in testing
\[
H_0 : C^{(h)}=C_1^{(h)}\widetilde\otimes C_2^{(h)}
\quad {\rm vs.}
 \quad H_A : C^{(h)}\ne C_1^{(h)}\widetilde\otimes C_2^{(h)}.
\]
In the following theorems notice that Theorem~\ref{thm:1} and
Theorem~\ref{thm:2} hold without the assumption of separability, i.e.
they hold under $H_0$ and under $H_A$. These two theorems are used to
establish the behavior of our test statistic under both the null,
Theorem~\ref{t:hatT},  and the alternative, Theorem
\ref{t:under:alt}.  Under the alternative both $C_1^{(h)}$ and
$C_2^{(h)}$ are still defined as partial traces of $C$, it is just
that their tensor product no longer recovers the original $C^{(h)}$.
\textcolor{black} {Before we state our theoretical results, we mention the asymptotic distribution of $\widehat{C}^{(h)}$,
which is the key to proof Theorem \ref{thm:1}.
It follows from Theorem 3 of \citetext{kokoszka:reimherr:2013}
that under  Assumption~\ref{a:main},
\begin{equation*} \label{e:KR:1}
\sqrt{N}(\widehat{C}^{(h)} - C^{(h)})
\xrightarrow{\mathcal{L}} N ( 0, \boldsymbol{\Gamma}^{(h)}),
\end{equation*}
 where $\boldsymbol{\Gamma}^{(h)}$ is given by
\begin{equation} \label{e:Gamma:1}
\boldsymbol{\Gamma}^{(h)} = \bR_0^{(h)} + \sum_{i=1}^\infty [\bR_i^{(h)} + (\bR_i^{(h)})^*]
\quad \text{with }
\bR_i^{(h)} =  \E [ ( \boldsymbol{X}_{1} \otimes \boldsymbol{X}_{1+h}- C^{(h)}) \otimes ( \boldsymbol{X}_{1+i} \otimes \boldsymbol{X}_{1+i+h} - {C}^{(h)})].
\end{equation}
Here $(\bR_i^{(h)})^*$ denotes the adjoint of $\bR_i^{(h)}$. Since we have the asymptotic distribution of $\widehat{C}^{(h)}$, and recalling that $\widehat{C}_1^{(h)}$ and $\widehat{C}_2^{(h)}$ are functions of $\widehat{C}^{(h)}$ from equation \eqref{e:hat-C12}, we can use the Delta method to prove the following theorem, details of the proof of Theorem \ref{thm:1} are given in Section~\ref{s:p}
of  Supporting Information.
}

\begin{theorem}\label{thm:1}
Under Assumption \ref{a:main}, one can explicitly define
a long--run covariance operator $\boldsymbol{W}^{(h)}$ such that
\begin{align*}
\sqrt{N} \left( \begin{matrix}
\widehat{C}_1^{(h)}- C_1^{(h)} \\
\widehat{C}_2^{(h)}- C_2^{(h)} \\
\widehat{C}^{(h)} - C^{(h)}
\end{matrix}
\right)
\xrightarrow{\mathcal{L}} N( {\bf 0}, \boldsymbol{W}^{(h)}).
\end{align*}
The definition of $\boldsymbol{W}^{(h)}$  is given in
formula \eqref{e:wmatrix} of  Supporting Information.
\end{theorem}

Armed with Theorem~\ref{thm:1},
we can derive the asymptotic distribution of
$\widehat{C}_1^{(h)}\widetilde\otimes\widehat{C}_2^{(h)}-\widehat{C}^{(h)}$.

\begin{theorem}\label{thm:2}
Under Assumption \ref{a:main},
\begin{equation*}
\sqrt{N} ( (\widehat{C}_1^{(h)}\widetilde \otimes \widehat{C}_2^{(h)}
- \widehat{C}^{(h)}) - (C_1^{(h)}\widetilde \otimes C_2^{(h)}  - C^{(h)}))\xrightarrow{\mathcal{L}}
N ( 0, \boldsymbol{Q}^{(h)})
\end{equation*}
The covariance operator  $\boldsymbol{Q}^{(h)} \in \cS(\cS(\mbH_1
\otimes \mbH_2))$ is defined in formula  \eqref{e:Qmatrix}
of  Supporting Information.
\end{theorem}

As a corollary, we obtain  the asymptotic distribution of
$\widehat{C}_1^{(h)}\widetilde\otimes\widehat{C}_2^{(h)}-\widehat{C}^{(h)}$ under $H_0$.
\begin{corollary}\label{corollary:1}
Suppose Assumption~\ref{a:main} holds. Then, under $H_0$,
\begin{equation*}
\sqrt{N}  (\widehat{C}_1^{(h)}\widetilde \otimes \widehat{C}_2^{(h)}
- \widehat{C}^{(h)})\xrightarrow{\mathcal{L}}
N ( 0, \boldsymbol{Q}^{(h)}),
\end{equation*}
\textcolor{black}{where the covariance operator $\boldsymbol{Q}^{(h)}
$ is the same as in Theorem~\ref{thm:2}.}
\end{corollary}

As noted above, in the context of \eqref{e:panel}, $\bQ^{(h)} \in
\cS(\cS(L_2^S))$, i.e. it is a Hilbert-Schmidt operator acting on a
space of Hilbert-Schmidt operators over $L_2^S$.  The following result
is a direct consequence of Theorem~\ref{thm:2}.  While the weighted
chi--square expansion is standard, to compute the weights, the
operator $\bQ^{(h)}$ must be estimated, so $\boldsymbol{W}^{(h)}$ must
be estimated. Formula \eqref{e:wmatrix} defining
$\boldsymbol{W}^{(h)}$ is new and nontrivial.

\begin{theorem} \label{t:hatT} Suppose Assumption~\ref{a:main}
 holds. Let $\bQ^{(h)}$ be the covariance
  operator appearing in Theorem \ref{thm:2}, whose  eigenvalues are
  $\gamma_1,\gamma_2,\dots$.  Then, under $H_0$, as $N\to \infty$,
\[
\widehat T \xrightarrow{\mathcal{L}}
\sum_{r=1}^\infty \gamma_r Z_r^2,
\]
where the $Z_r$ are iid standard normal.
\end{theorem}

To describe the behavior of the test statistic under the alternative,
some specific  form of the alternative must be assumed, as the
violation of \eqref{e:H0-C} can take many forms. A  natural
approach corresponding to a fixed alternative to
$C_1^{(h)}\widetilde \otimes C_2^{(h)}  - C^{(h)} = 0$, is to assume
that
\begin{equation} \label{e:HA}
C_1^{(h)}\widetilde \otimes C_2^{(h)}  - C^{(h)} = : \Delta \neq 0.
\end{equation}

\begin{theorem}\label{t:under:alt}
Suppose Assumption \ref{a:main} holds. If \eqref{e:HA} holds, then
\[
\widehat{T}=N\|\Delta\|^2 + O_P(N^{1/2})
\xrightarrow{\mathcal{P}} \infty.
\]
\end{theorem}

\medskip

In our applications, $\boldsymbol{X}_n \in \mbH_1 \otimes \mbH_2$,
where $\mbH_1=\mbR^S$ and $\mbH_2=L^2([0,1])$.  Therefore, in
practice, we must first project these random elements onto a truncated
basis by using a dimension reduction procedure. Note that
$\mbH_1=\mbR^S$ is already finite. However, if the number of
coordinates in the panel is large, then a dimension reduction in
$\mbH_1=\mbR^S$ is also recommended.  Here we present the general case
where we use dimension reduction in both $\mbH_1=\mbR^S$ and
$\mbH_2=L^2([0,1])$.  \textcolor{black} {The truncated basis is of the form
$\hat{u}_k\otimes \hat{v}_j$ with $1 \le k \le K,\ 1 \le j \le J$
where $K<S$ and $J < \infty$. In our implementation,  the
$\hat{u}_k$ and the  $\hat{v}_j$ are the empirical principal components.
We can approximate each
$\boldsymbol{X}_n \in \mbH_1 \otimes \mbH_2$ by a $K\times J$ random
matrix $\bZ_n \in \mbR^{K\times J}$, where $\bZ_n(k,j)=\langle
\boldsymbol{X}_n,\hat{u}_k\otimes \hat{v}_j\rangle$, $1 \le k \le K,\
1 \le j \le J$.  Therefore, from now on, we work with observations in
the form of random $K\times J$ matrices defined as
\[
\bZ_n = [z_{kj;n},\  1 \le k \le K,\ 1 \le j \le J],
\]
where
$z_{kj;n}=\langle \boldsymbol{X}_n,
\hat{u}_k\otimes \hat{v}_j\rangle$.
}%textcolor
Let $\widehat{T}_F$ be the truncated test statistic $\widehat{T}$, i.e.
\[
\widehat T_F
= N \| \widehat C^{(h)}_{1,K} \widetilde\otimes \widehat C^{(h)}_{2,J}
- \widehat C^{(h)}_{KJ} \|_F^2,
\]
where $\widehat C^{(h)}_{1,K}$ is a $K\times K$ matrix, $\widehat
C^{(h)}_{2,J}$ is a $J\times J$ matrix, $\widehat C^{(h)}_{KJ}$ is
a fourth order array of dimension $K\times J\times K\times J$, and
$\|\cdot\|_F$ is the Frobenius norm, which is the Hilbert--Schmidt
norm in finite dimensions.  Finally, let $\bQ^{(h)}_{KJ}$ be the
truncated covariance operator $\bQ^{(h)}$, i.e. $\bQ^{(h)}_{KJ}$ is
the asymptotic covariance operator in the convergence
\[
\sqrt{N} ( (\widehat{C}^{(h)}_{1,K}\widetilde \otimes \widehat{C}^{(h)}_{2,J}  - \widehat{C}^{(h)}_{KJ}) - (C^{(h)}_{1,K}\widetilde \otimes C^{(h)}_{2,J}  - C^{(h)}_{KJ}))\xrightarrow{\mathcal{L}}
N ( 0, \boldsymbol{Q}^{(h)}_{KJ}).
\]
Note that $\bQ^{(h)}_{KJ}$ is an  array of order eight with finite dimensions,
$\bQ^{(h)}_{KJ} \in \mbR^{K\times J\times K\times J\times K\times
J\times K\times J}$. More details are given in Remark~\ref{r:in:practice}
in  Supporting Information.
As a finite array, it has only a finite number
of eigenvalues, which with
denote $\gamma_1^\dagger,\gamma_2^\dagger,\dots,\gamma_R^\dagger$.
The arguments leading to Theorem~\ref{t:hatT}  show that
under  $H_0$, as $N\to \infty$,
\begin{equation} \label{e:dist}
\widehat {T}_F \xrightarrow{\mathcal{L}}
\sum_{r=1}^R \gamma_r^\dagger Z_r^2,
\end{equation}
where the $Z_r$ are iid standard normal. The asymptotic argument
needed to establish \eqref{e:dist} relies on the bounds
$\| \hat u_k - u_k \| = O_P(N^{-1/2})$ and
$\| \hat v_j  - v_j  \| = O_P(N^{-1/2})$, which hold under
Assumption~\ref{a:main}. It is  similar to the technique
used in the proof of Theorem 4 in  \citetext{constantinou:2017},
so it is omitted.

\section{Details of implementation}\label{s:det}
Recall that we assume that all functions have been rescaled so that
their domain is the unit interval $[0,1]$, and that they have mean
zero. The testing procedure consists of dimension reduction in time
and, for large panels, a further dimension reduction in
coordinates. After reducing the dimension our "observations" are of
the form of $K\times J$ matrices which are used to compute the estimators
we need to perform our test.  The remainder of this section
explains the details in an algorithmic form.  The reader will notice
that most steps have obvious variants, for example, different weights
and bandwidths can be used in Step 6.  Procedure~\ref{proc:both-PC}
describes the exact implementation used in Sections \ref{s:sim} and
\ref{s:appl}.

\begin{procedure} \label{proc:both-PC} \hspace*{15cm}\\

\noindent {\bf 1.}
[{\it Pool across $s$ to get estimated temporal FPCs.}]  Under the
assumption of separability, i.e., under the $H_0$ stated in
Section~\ref{s:ap}, the optimal functions used for temporal dimension
reduction are the same for each member (coordinate) of the panel;
information can then be pooled across the coordinates to get better
estimates of these functions. In other words, under separability,  we can
use simultaneously all the $N\times S$ functions to compute the  temporal
FPCs $\hat{v}_1,\dots,\hat{v}_J$ as  the eigenfunctions
of the covariance function
\[
\hat {c}_2(t,t^\prime)
=\frac{1}{NS}\sum_{n=1}^{N}\sum_{s=1}^SX_{ns}(t)X_{ns}(t^\prime).
\]

\noindent {\bf 2.} Approximate each curve  $X_{ns}(t)$ by
\[
X^{(J)}_{ns}(t) = \sum_{j=1}^J
\xi_{nsj} \hat{v}_j(t),
\]
where $\xi_{nsj}=\langle X_{ns}(t), \hat{v}_j(t)\rangle$.
Construct $S\times J$ matrices $\bXi_n$ defined as
\begin{equation*} \label{e:Xi-n}
\bXi_n = [\xi_{nsj},\  1 \le s \le S,\ 1 \le j \le J],
\end{equation*}
\noindent where $J$ is chosen large enough so that the first $J$ FPCs
explain at least $85\%$ of the variance.
This is Functional Principal Components Analysis carried out on the
pooled (across coordinates) sample.\\

\noindent {\bf 3.} [{\it Pool across time to get panel PCs.}]
Under the assumption of separability the panel principal components
are the same for each time. In other words the panel PCs are the
principal components of the following covariance matrix:
\[
\hat {c}_1(s,s')
= \frac{ \sum_n^{N} \int X_{ns}(t) X_{ns'}(t) \ dt}{N \tr( \widehat {C} )}.
\]
However, since we have already reduced the dimension of the observed
functions,  the panel PCs
$\hat{{u}}_1,\dots,\hat{{u}}_K$ are the
principal components of the covariance matrix
\[
\tilde c_{1}(s , s') = \frac{1}{NJ}\sum_{n=1}^{N}
 \sum_{j=1}^J \frac{\xi_{nsj} \xi_{ns' j}}{\lambda_{j}}.
\]
\\
\noindent {\bf 4.} Approximate each row
${\xi}_{n\cdot j}=(\xi_{n1j}, \xi_{n2j}, \dots, \xi_{nSj})$
of the $\bXi_n$ matrices  by
\[
{\xi}^{(K)}_{n\cdot j} = \sum_{k=1}^K z_{kj;n}\hat{{u}}_k, \ \ \
z_{kj;n}=\langle {\xi}_{n\cdot j},\hat{{u}}_k\rangle.
\]
Construct the $K\times J$ matrices
$\bZ_n = [z_{kj;n},\ 1 \le k \le K,\ 1 \le j \le J]$,  where $K$ is
chosen large enough so that the first $K$ eigenvalues explain at least
$85\%$ of the variance.  This is a multivariate PCA on the pooled
(across time) variance adjusted sample.\\

\noindent If the number of panel coordinates is small, then a
multivariate dimension reduction is not necessary, so one can skip
steps $3$ and $4$ and use the $\bXi_n$ matrices instead of the $\bZ_n$
matrices, and replace $K$ with $S$ in the following steps.
\textcolor{black}{The dimension reduction steps
reduce the computational time and the memory requirements
by reducing the matrix size  the 4D and 8D covariance tensors. } \\

\noindent {\bf 5.} Approximate covariance \eqref{e:hat-Ch} by
the  fourth order array of dimensions
$K\times J\times K \times J$
\[
\widehat{C}^{(h)}_{KJ}=\frac{1}{N-h}
\sum_{n=1}^{N-h} \boldsymbol{Z}_n\otimes \boldsymbol{Z}_{n+h}.
\]
Approximate  $\widehat C_1^{(h)}$ and $\widehat C_2^{(h)}$
in \eqref{e:hat-C12} by
\[
\widehat{C}^{(h)}_{1,K}(k,k^\prime) = \frac{\sum_{j=1}^J\widehat{C}^{(h)}_{KJ}(k,j,k^\prime,j)}{\sum_{k=1}^K \sum_{j=1}^J C^{(h)}_{KJ}(k,j,k,j)}
 \quad  \text{and}
\quad \widehat{C}^{(h)}_{2,J}(j,j^\prime) = \sum_{k=1}^K\widehat{C}^{(h)}_{KJ}(k,j,k,j^\prime),
\]
where $\widehat{C}^{(h)}_{1,K}$ is a $K\times K$ matrix and
$\widehat{C}^{(h)}_{2,J}$ is a $J\times J$ matrix.\\

\noindent {\bf 6.} \textcolor{black} {Calculate the estimators $\widehat{\bR}^{(h)}_{0,KJ}, \widehat{\bR}^{(h)}_{i,KJ}, (\widehat{\bR}^{(h)}_{i,KJ})^* \in \mbR^{K\times J\times K\times J\times K\times J\times K\times J}$, by using
\begin{align}\label{a:R}
\begin{split}
 & \widehat{\bR}^{(h)}_{0,KJ}=\frac{1}{N-h} \sum_{n=1}^{N-h} [ ( \boldsymbol{Z}_n\otimes \boldsymbol{Z}_{n+h} - \widehat{C}^{(h)}_{KJ}) \otimes ( \boldsymbol{Z}_n\otimes \boldsymbol{Z}_{n+h} - \widehat{C}^{(h)}_{KJ})], \\
&\widehat{\bR}^{(h)}_{i,KJ}=\frac{1}{N-i-h} \sum_{n=1}^{N-i-h} [ ( \boldsymbol{Z}_n\otimes \boldsymbol{Z}_{n+h} - \widehat{C}^{(h)}_{KJ}) \otimes ( \boldsymbol{Z}_{n+i}\otimes \boldsymbol{Z}_{n+i+h} - \widehat{C}^{(h)}_{KJ})],\\
 & (\widehat{\bR}^{(h)}_{i,KJ})^*= \frac{1}{N-i-h} \sum_{n=1}^{N-i-h} [ ( \boldsymbol{Z}_{n+i}\otimes \boldsymbol{Z}_{n+i+h} - \widehat{C}^{(h)}_{KJ}) \otimes ( \boldsymbol{Z}_n\otimes \boldsymbol{Z}_{n+h} - \widehat{C}^{(h)}_{KJ})].
\end{split}
\end{align}}

\noindent {\bf 7.} \textcolor{black}{Calculate the estimator $\widehat{\boldsymbol{\Gamma}}^{(h)}_{KJ} \in \mbR^{K\times J\times K\times J\times K\times J\times K\times J}$, by using the following Bartlett-type estimator:
\begin{equation} \label{e:Bartlett:Gamma}
\widehat{\boldsymbol{\Gamma}}^{(h)}_{KJ} = \widehat{\bR}^{(h)}_{0,KJ} + \sum_{i=1}^{N-h-1} \omega_i (\widehat{\bR}^{(h)}_{i,KJ} + (\widehat{\bR}^{(h)}_{i,KJ})^* ),
\end{equation}
where $\widehat{\bR}^{(h)}_{0,KJ}, \widehat{\bR}^{(h)}_{i,KJ}, (\widehat{\bR}^{(h)}_{i,KJ})^*$ are defined in equation \eqref{a:R} and
the  $\omega_i$ are the Bartlett's weights, i.e.,
\[
\omega_i =
\begin{cases}
1-\frac{i}{1+q}, & \mbox{if} \hspace{0.2cm} i \le q\\
0, & \mbox{otherwise},
\end{cases}
\]
with $i$ being the number of lags and $q$ is the bandwidth which is
assumed to be a function of the sample size, i.e., $q=q(N)$. In our
simulations, in Section ~\ref{s:sim}, we use the formula
$q\approx1.1447(\frac{N}{4})^{1/3}$ (\citetext[Chapter 16]{HKbook}).}\\

\textcolor{black}{Note that the estimators $\widehat{\bR}^{(h)}_{0,KJ}, \widehat{\bR}^{(h)}_{i,KJ}, (\widehat{\bR}^{(h)}_{i,KJ})^*,\widehat{\boldsymbol{\Gamma}}^{(h)}_{KJ}$ defined in steps 6 and 7 are the truncated analogs of the estimators $\widehat{\bR}^{(h)}_{0}, \widehat{\bR}^{(h)}_{i}, (\widehat{\bR}^{(h)}_{i})^*,\widehat{\boldsymbol{\Gamma}}^{(h)}$, which can be obtained by simply changing $\boldsymbol{Z_n}$ with $ \boldsymbol{X_n}$ in equation \eqref{a:R}.\\
}

\noindent {\bf 8.}
Estimate the arrays $\bW^{(h)}_{KJ}$ (the truncated analog
of $\bW^{(h)}$) and $\bQ^{(h)}_{KJ}$ defined in Section~\ref{s:da}.
Details are given in Remark~\ref{r:in:practice} in
Supporting Information.\\

\noindent {\bf 9.} Calculate the P--value using the limit distribution
specified in \eqref{e:dist}.\\

\end{procedure}

Step 2 can be easily implemented using {\tt R} function {\tt pca.fd},
and step 3 by using {\tt R} function {\tt prcomp}. The matrix
$\bQ^{(h)}_{KJ}$ can be computed using the {\tt R} package {\tt
  tensorA} by \citetext{boogaart:gerald:2007}.

\section{A simulation study}\label{s:sim}
The purpose of this section is to provide information on the
performance of our  test procedure in finite samples. We first
 comment on the performance
of existing tests.
\citetext{constantinou:2017} derived several separability  tests based on
the assumption of independent $\boldsymbol{X}_n$.  For the functional
panels which exhibit temporal dependence (we define them below), the
empirical sizes are close to zero; the tests of
\citetext{constantinou:2017} are too conservative to be usable, unless
we have  independent replications of the spatio--temporal
structure. \citetext{aston2017tests} proposed three tests,
also for independent $\boldsymbol{X}_n$.  In the presence of temporal
dependence, their tests are not useable either; they can severely
overreject, the empirical size  can approach 50\% at the
nominal level of 5\%. We give some specific numbers at the end of this
section.

For our empirical study,
we simulate functional panels as the moving average process
\[
X_{ns}(t)= \sum_{s'=1}^S  \Psi_{ss^\prime}
[e_{ns'}(t)+e_{n-1s'}(t)],
\]
which is a 1-dependent functional time series.  Direct verification,
shows that it is separable as long the $e_{ns}(t)$ are separable.
We generate $e_{ns}(t)$ as Gaussian processes with the
following covariance function, which is a modified version of Example 2 of
\citetext{cressie:huang:1999}:
\begin{equation} \label{e:cov1}
\sigma_{ss^\prime}(t,\tp)
=\frac{\sigma^2}{(a|t-t^\prime|+1)^{1/2}}
\exp\left(-\frac{b^2[|s - s^\prime|/(S-1)]^{2}}{(a|t - \tp|+1)^c}\right).
\end{equation}
In this covariance function, $a$ and $b$ are nonnegative scaling
parameters of time and space, respectively, and $\sigma^2>0$ is an
overall scale parameter. The most important parameter is the
separability parameter $c$ which takes values in $[0,1]$. If $c=0$,
the covariance function is separable, otherwise it is not. We set
$a=3$, $b=2$, $\sigma^2=1$. To simulate the functions,
we use $T=50$ time points equally spaced on $[0,1]$,
and  $S \in
\{4,6,8,10,12,14\}$ coordinates in the panel.
The MA coefficients  are taken as:
\[
\Psi_{ss^\prime}
=\exp\left(-\frac{25(s - s^\prime)^2}{(S-1)^{2}}\right).
\]
Notice that in the covariance
above, the differences in the coordinates of the panel, i.e. $|s -
s^\prime|$, are rescaled to be within the interval $[0,1]$, i.e. we use
$|s - s^\prime|/(S-1)$.

We set
\[
c=0 \ {\rm under} \ H_0;
\ \ \ c=1  \ {\rm under} \ H_A.
\]

We  consider two different cases. The first one
with dimension reduction only in time and the second one with
dimension reduction in both time and coordinates. For each case,  we
study two different scenarios. The first scenario is under the null
hypothesis (separability) and the second scenario under the
alternative hypothesis. We consider different
numbers of temporal FPCs, $J$, in the first case and different
numbers of coordinate PCs, $K$, and temporal FPCs, $J$, in the
second case.  We will also consider different values for the series length
$N$. All empirical rejection rates are based on
one thousand replications, so their SD is about 0.7 percent for size
(we use the nominal significance level of 5\%), and about two percent for
power.

\subsection{Case 1: dimension reduction in time only}\label{ss:dr1}
We examine the effect of the series length $N$ and the
number of principal components $J$ on the empirical
size (Table~\ref{t:n})  and power (Table~\ref{t:al})
for $S \in \{4,6,8\}$. Each table reports the
rejection rates in percent. In parentheses,   the
proportion of variance explained by the $J$ PCs is given.

\begin{table*}[ht]\centering
\begin{tabular}{@{}cccccccccccc@{}}\toprule
&\multicolumn{3}{c}{$N=100$} & & \multicolumn{3}{c}{$N=150$} & & \multicolumn{3}{c}{$N=200$} \\
\cmidrule{2-4} \cmidrule{6-8} \cmidrule{10-12}
&$J=2$ & $J=3$ & $J=4$ && $J=2$ & $J=3$ & $J=4$ && $J=2$ & $J=3$ & $J=4$ \\ \midrule

$S=4$ & 5.5 &6.4 & 5.0 && 5.9& 5.7 & 5.3 && 6.5 & 5.1 & 5.5 \\
            & (87\%)&(90\%)&(94\%)&&(85\%)&(90\%)&(92\%)&&(87\%)&(90\%)&(92\%)\\
$S=6$ & 5.6 & 5.9 & 5.3 && 6.2 & 5.3 & 4.7 && 5.6 & 6.2 & 5.1 \\
            & (85\%)&(91\%)&(93\%)&&(85\%)&(91\%)&(93\%)&&(86\%)&(91\%)&(92\%)\\
$S=8$ & 5.4 &6.0 & 7.5 && 4.8 & 5.8 & 6.6 && 6.0 & 5.7 & 6.1 \\
            & (87\%)&(89\%)&(94\%)&&(86\%)&(91\%)&(94\%)&&(85\%)&(89\%)&(93\%)\\

\bottomrule
\end{tabular}
\vspace{0.2cm}
\caption{Rejection rates under $H_0$ ($c=0$) at the nominal
5 percent level;  $J$ is the number of temporal PCs.
 The explained variance of the temporal FPCA
is given in parentheses. \label{t:n}}
\end{table*}

\begin{table*}[ht]\centering
\begin{tabular}{@{}cccccccccccc@{}}\toprule
&\multicolumn{3}{c}{$N=100$} & & \multicolumn{3}{c}{$N=150$} & & \multicolumn{3}{c}{$N=200$} \\
\cmidrule{2-4} \cmidrule{6-8} \cmidrule{10-12}
&$J=2$ & $J=3$ & $J=4$ && $J=2$ & $J=3$ & $J=4$ && $J=2$ & $J=3$ & $J=4$ \\ \midrule

$S=4$ & 67.6 &90.6 & 95.1 && 91.9& 99.3 & 99.8 && 98.2 & 100 & 100 \\
            & (86\%)&(90\%)&(94\%)&&(87\%)&(90\%)&(93\%)&&(87\%)&(92\%)&(94\%)\\
$S=6$ & 54.5& 79.7 & 89.0 && 80.7 & 97.9 & 99.3 && 94.5 & 99.7 & 100 \\
            & (88\%)&(91\%)&(94\%)&&(85\%)&(91\%)&(94\%)&&(88\%)&(92\%)&(94\%)\\
$S=8$ & 45.2 &74.9 & 85.2 && 75.1 & 96.8 & 98.7 && 91.5 & 99.9 & 100 \\
            & (89\%)&(91\%)&(94\%)&&(89\%)&(92\%)&(94\%)&&(88\%)&(92\%)&(94\%)\\
\bottomrule
\end{tabular}
\vspace{0.2cm}
\caption{Empirical power ($c=1$);
$J$ is  the number of temporal PCs.  The explained variance of the
temporal FPCA is given in parentheses.  \label{t:al}}
\end{table*}

In Table \ref{t:n}, we can see that the size of  our test  is robust to
the number of the principal components used. This is a very desirable
property, as in all procedures of FDA there is some uncertainty about  the
optimal number of FPCs that should be used.
While still within two standard errors of the nominal size,
the empirical size becomes inflated for $S=8$. We recommend dimension
reduction in panel coordinates if $S\ge 10$.
In Table \ref{t:al}, we  see that the empirical power increases
as $N$ and $J$
increase.  The power increase with $N$ is  expected;
its increase with $J$ reflects the fact that projections on  larger
subspaces better capture a  departure from $H_0$. However,
$J$ cannot be chosen too large so as not to increase the dimensionality
of the problem, which negatively affects the empirical size.

\subsection{Case 2: dimension reduction in both time and
panel coordinates}\label{ss:dr2} The general setting is
the same as in Section~\ref{ss:dr1}, but we consider larger
panels, $S \in \{10,12,14\}$, and reduce their
dimension  to $K \in \{2,3,4\}$ coordinates.
The proportion of variance  explained is now  computed as
\begin{equation} \label{e:cpv}
{\rm CPV}(J, K) =
\frac{\sum_{j=1}^J
 \lambda_j}{\sum_{j=1}^\infty \lambda_j} \times \frac{\sum_{k=1}^K
 \mu_k}{\sum_{k=1}^S \mu_k},
\end{equation}
where the $\lambda_1, \lambda_2, \dots$, and $\mu_1,
 \mu_2, \dots, \mu_S$ are, respectively,
the estimated  eigenvalues of the time and panel PCA's.

\begin{table*}[ht]\centering
{\footnotesize
\begin{tabular}{@{}c|cccccccccccc@{}}\toprule
 \multicolumn{1}{c}{}&\multicolumn{1}{c}{}&\multicolumn{3}{c}{$N=100$} & & \multicolumn{3}{c}{$N=150$} & & \multicolumn{3}{c}{$N=200$} \\
\cmidrule{3-5} \cmidrule{7-9} \cmidrule{11-13}
\multicolumn{1}{c}{}&&$J=2$ & $J=3$ & $J=4$ && $J=2$ & $J=3$ & $J=4$ && $J=2$ & $J=3$ & $J=4$ \\ \midrule

&$K=2$ & 6.4 &6.2 & 6.1 && 6.1& 5.2 & 4.3 && 5.8 & 5.5 & 5.6 \\
              & & (80\%)&(84\%)&(90\%)&&(80\%)&(85\%)&(88\%)&&(80\%)&(84\%)&(88\%)\\
$S=10$&$K=3$ & 6.1 & 4.8 & 5.6 && 5.0 & 5.5 & 4.7 && 5.3 & 6.1 & 5.0 \\
              & & (84\%)&(88\%)&(94\%)&&(83\%)&(89\%)&(92\%)&&(85\%)&(89\%)&(92\%)\\
&$K=4$ & 5.9 &6.3 & 4.7 && 5.2 & 5.8 & 5.3 && 6.1 & 6.1 & 5.8 \\
              && (84\%)&(90\%)&(92\%)&&(85\%)&(90\%)&(92\%)&&(84\%)&(90\%)&(92\%)\\
\hline
&$K=2$ & 6.3 &6.4 & 6.1 && 5.6& 6.0 & 6.2 && 4.6 & 6.3 & 6.4 \\
              & & (83\%)&(88\%)&(90\%)&&(83\%)&(86\%)&(89\%)&&(80\%)&(87\%)&(88\%)\\
$S=12$&$K=3$ & 6.1 & 5.9 & 5.1 && 5.0 & 5.6 & 6.0 && 6.1 & 4.8 & 6.1 \\
              & & (87\%)&(91\%)&(93\%)&&(87\%)&(90\%)&(92\%)&&(85\%)&(90\%)&(93\%)\\
&$K=4$ & 6.0 &5.4 & 5.0 && 6.3 & 6.1 & 6.0 && 5.0 & 6.5 & 5.7 \\
              && (87\%)&(91\%)&(93\%)&&(85\%)&(90\%)&(93\%)&&(86\%)&(90\%)&(93\%)\\
\hline
&$K=2$ & 6.4 &5.2 & 4.5 && 6.2& 5.8 & 5.0 && 5.6 & 6.6 & 5.3 \\
              & & (82\%)&(87\%)&(89\%)&&(82\%)&(87\%)&(88\%)&&(82\%)&(86\%)&(89\%)\\
$S=14$&$K=3$ & 6.0 & 5.2 & 4.7 && 4.4 & 6.2 & 6.0 && 4.2 & 5.6 & 6.2 \\
              & & (85\%)&(90\%)&(92\%)&&(83\%)&(88\%)&(92\%)&&(84\%)&(90\%)&(93\%)\\
&$K=4$ & 6.0 &5.1 & 4.6 && 6.3 & 5.5 & 5.7 && 6.5 & 5.7 & 5.6 \\
              && (85\%)&(90\%)&(93\%)&&(86\%)&(90\%)&(93\%)&&(87\%)&(89\%)&(91\%)\\
\bottomrule
\end{tabular}
\vspace{0.3cm}
\caption{Rejection rates under $H_0$ ($c=0$);
 $K$ is the reduced panel dimension and $J$ the number of temporal PCs.
The explained variance of the dimension reduction is given in parentheses.
\label{t:n1} }
}%footnotesize
\end{table*}

\begin{table*}[ht]\centering
{\footnotesize
\begin{tabular}{@{}c|cccccccccccc@{}}\toprule
 \multicolumn{1}{c}{}&\multicolumn{1}{c}{}&\multicolumn{3}{c}{$N=100$} & & \multicolumn{3}{c}{$N=150$} & & \multicolumn{3}{c}{$N=200$} \\
\cmidrule{3-5} \cmidrule{7-9} \cmidrule{11-13}
\multicolumn{1}{c}{}&&$J=2$ & $J=3$ & $J=4$ && $J=2$ & $J=3$ & $J=4$ && $J=2$ & $J=3$ & $J=4$ \\ \midrule

&$K=2$ & 31.3 &49.2 & 60.8 && 50.5 & 78.0 & 85.9 && 67.3 & 92.4 & 96.4 \\
              & & (85\%)&(85\%)&(83\%)&&(81\%)&(85\%)&(84\%)&&(81\%)&(82\%)&(84\%)\\
$S=10$&$K=3$ & 43.0 & 73.3 & 81.7 && 72.3 & 95.2 & 98.2 && 90.6 & 99.5 & 100 \\
              & & (87\%)&(90\%)&(93\%)&&(85\%)&(91\%)&(93\%)&&(86\%)&(91\%)&(93\%)\\
&$K=4$ & 42.2 & 73.5 & 84.7 && 71.1 & 96.6 & 98.5 && 90.6 & 99.8 & 100 \\
              && (85\%)&(92\%)&(93\%)&&(87\%)&(92\%)&(94\%)&&(88\%)&(90\%)&(92\%)\\
\hline
&$K=2$ & 30.8 & 49.7 & 57.4 && 46.6 & 76.5 & 87.2 && 67.7 & 91.8 & 95.7 \\
              & & (82\%)&(83\%)&(85\%)&&(81\%)&(83\%)&(84\%)&&(82\%)&(83\%)&(86\%)\\
$S=12$&$K=3$ & 42.9 & 72.5 & 82.8 && 67.1 & 94.8 & 98.8 && 89.5 & 99.5 & 99.9 \\
              & & (89\%)&(91\%)&(94\%)&&(88\%)&(92\%)&(93\%)&&(87\%)&(93\%)&(93\%)\\
&$K=4$ & 43.3 &72.0 & 82.9 && 71.1 & 95.9 & 97.8 && 89.0 & 99.6 & 100 \\
              && (87\%)&(92\%)&(94\%)&&(86\%)&(92\%)&(94\%)&&(86\%)&(91\%)&(93\%)\\
\hline
&$K=2$ & 27.7 & 46.2 & 55.0 && 47.7 & 74.9 & 82.8 && 69.0 & 90.6 & 94.0 \\
              & & (86\%)&(84\%)&(84\%)&&(82\%)&(83\%)&(84\%)&&(81\%)&(84\%)&(87\%)\\
$S=14$&$K=3$ & 39.2 & 66.6 & 81.3 && 67.5 & 91.0 & 93.4 && 88.1 & 94.4 & 94.1 \\
              & & (87\%)&(92\%)&(93\%)&&(89\%)&(91\%)&(93\%)&&(88\%)&(90\%)&(93\%)\\
&$K=4$ & 43.7 &70.4 & 78.9 && 70.5 & 91.1 & 93.7 && 88.2 & 94.4 & 95.7 \\
              && (87\%)&(92\%)&(94\%)&&(88\%)&(91\%)&(94\%)&&(88\%)&(93\%)&(94\%)\\
\bottomrule
\end{tabular}
\vspace{0.2cm}
\caption{Empirical power ($c=1$); $K$ and $J$ are
as in Table~\ref{t:n1}.  The explained variance of the dimension
reduction is given in parentheses.\label{t:al1}}
}%footnotesize
\end{table*}

Tables \ref{t:n1} and \ref{t:al1} show that the reduction of
the panel dimension does not negatively affect the properties of the
tests. The conclusions are the same as in Section~\ref{ss:dr1}.
Either approach leads to a test with well controlled size, which
is does not depend on $J$ ($J, K$) as long the the proportion of
explained variance remains within the generally recommended
range of  85\%--95\%. If $J=2$ or $K=2$ are used, this requirement
is generally not met, resulting in a size distortion, which is
however acceptable and decreases with $N$.

As noted at the beginning of this section, the tests of
\citetext{constantinou:2017} are too conservative, they almost never
reject under the null for all scenarios considered in this section.  The tests of
\citetext{aston2017tests} reject too often under the
null. For example, in the settings considered in Table~\ref{t:n1}, the
rejection rates for their asymptotic test, Gaussian parametric
bootstrap test, and Gaussian parametric bootstrap test using
Hilbert--Schmidt distance, range between $19.0\%-49.4\%$,
$14.6\%-32.2\%$ and $38.1\%-44.9\%$, respectively.
By contrast, the test derived in this paper, in its both versions
and under all reasonable choices of  tuning parameters, has
precise empirical size at the standard 5\% nominal level and useful
power.

\textcolor{black}{In Section B of  Supporting Information, we
show the results of other simulations which study the effect of  different
covariance functions, the magnitude of the departure from $H_0$,
and the lag $h$. They do not modify the general conclusion that
the test is reasonably well calibrated and has useful power.}

\section{Applications to pollution and stock market data}\label{s:appl}
We begin by applying  our  method to air quality data
studied by \citetext{constantinou:2017} under the
assumption that the monthly curves are iid. These curves however
form a time series, so it is important to check if a test that accounts
for the temporal dependence leads to  the same or a different conclusion.

The Environmental Protection Agency (EPA)
collects massive amounts of air quality data which are available
through its website
\url{http://www3.epa.gov/airdata/ad_data_daily.html}.  The records
consist of data for $6$ common pollutants, collected by outdoor
monitors in hundreds of locations across the United States.  The
number and frequency of the observations varies greatly by location,
but some locations have as many as 3 decades worth of daily
measurements.  We focus on nitrogen dioxide, a common
pollutant emitted by combustion engines and power stations.

\begin{figure}[ht]
\begin{center}
\includegraphics[width=14cm,height=8.5cm]{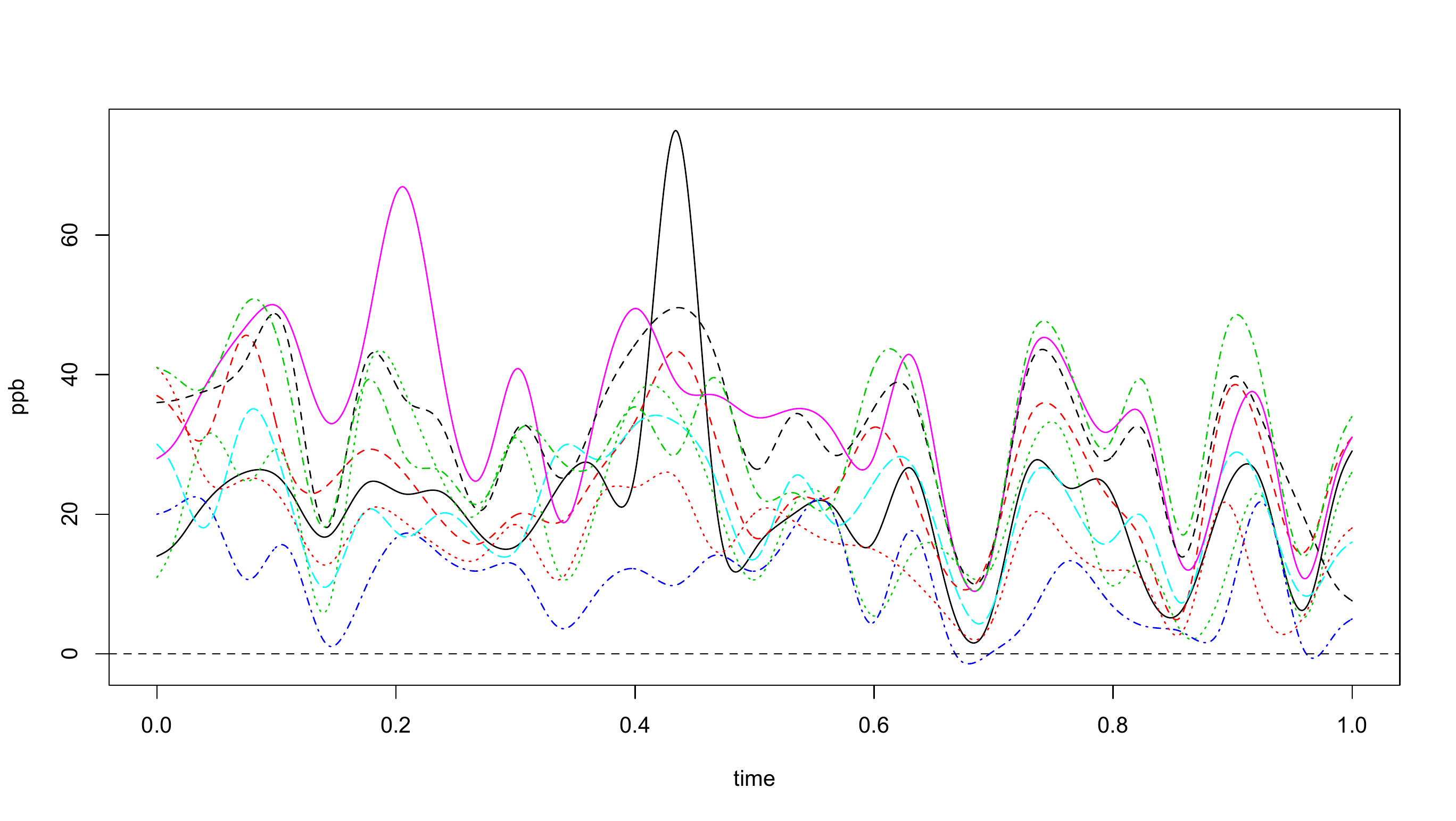}
  \caption{Maximum one--hour nitrogen dioxide curves for December
    2012 at the nine locations.
    \label{f:air_data}}
\end{center}
\end{figure}

We consider nine locations along the east coast that have relatively
complete records since 2000: Allentown, Baltimore, Boston, Harrisburg,
Lancaster, New York City, Philadelphia, Pittsburgh, and Washington D.C.
We use the data for the years 2000-2012 .  Each functional observation
$X_{ns}(t)$ consists of the daily maximum one-hour nitrogen dioxide
concentration measured in ppb (parts per billion) for day $t$, month
$n$ ($N=156$), and at location $s$. We thus have
a panel of $S=9$ functional time series
 (one at every location),  $X_{ns}(t), s=1,2,\dots, 9, n=1,2,\dots,156$.
Figure~\ref{f:air_data} shows the data for the nine locations for
December 2012.  Before the application of the test, the curves
were deseasonalized by  removing the monthly mean from each curve.

We applied both versions of Procedure \ref{proc:both-PC} (dimension in
time only and double dimension reduction).  Requiring 85\% to 95\% of
explained variance yielded the values $J, K=2,3, 4$, similarly as in
our simulated data example.  For all possible combinations, we
obtained P--values smaller than 10E-4.  This indicates a nonseparable
covariance function and confirms the conclusion obtained by
\citetext{constantinou:2017}; nonseparability is
an intrinsic feature of pollution data, simplifying the covariance
structure by assuming separability may lead to incorrect conclusions.

\medskip

\begin{figure}[ht]
\begin{center}
\includegraphics[width=14cm,height=8.5cm]{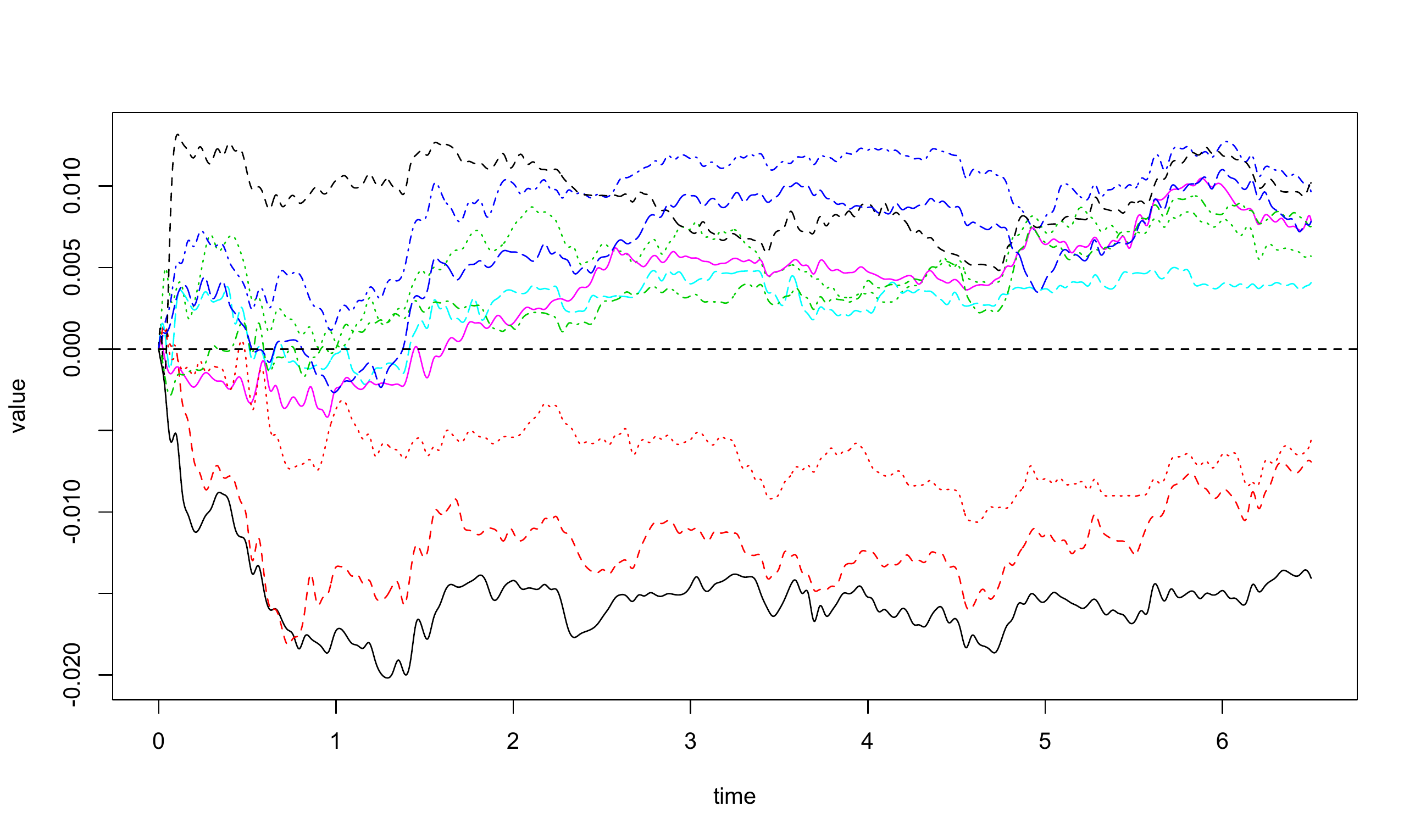}
  \caption{Cumulative intraday return
curves for the ten companies for April 2nd, 2007.
    \label{f:stock_data}}
\end{center}
\end{figure}

We now turn to an application to a stock portfolio.
Cumulative intradaily returns have recently been studied in several papers,
including \citetext{kokoszka:reimherr:2013pred},
\citetext{kokoszka:miao:zhang:2015} and \citetext{lucca:moench:2015}.
If $P_n(t)$ is the price of a stock at minute $t$ of the trading day $n$,
then the cumulative intraday return curve on day $n$ is defined
by
\[
R_n(t)=\log(P_n(t)) - \log(P_n(0)),
\]
where time $0$ corresponds to the opening of the market
(9:30 EST for the NYSE).
\citetext{horvath:kokoszka:rice:2014}
\textcolor{black} {did not find evidence against temporal stationarity
of such time series. The work of
\citetext{kokoszka:reimherr:2013pred} shows that
cumulative intradaily returns do not form an iid sequence.
(This can be readily verified by computing the ACF of squared scores.)
}%textcolor
 Figure~\ref{f:stock_data} shows the
curves $R_n$  for  ten companies on  April 2nd.
2007. This  portfolio of $S=10$ stocks produces a panel of functional
time series studied in this paper. We selected ten US blue chip
companies, and want to determine if the resulting panel can
be assumed to have a separable covariance function. The answer
is yes, as we now explain.

We consider stock values, recorded every minute, from October 10, 2001
to April 2, 2007 (1,378 trading days) for the following 10 companies:
Bank of America (BOA), Citi Bank, Coca Cola, Chevron Corporation
(CVX), Walt Disney Company (DIS), International Business Machines
(IBM), McDonald's Corporation (MCD), Microsoft Corporation (MSFT),
Walmart Stores (WMT) and Exxon Mobil Corporation Common (XOM). On each
trading day, there are 390 discrete observations.  There is an outlier
on August 26, 2004 for Bank of America, which is due to a stock
split. That day is discarded from further analysis, so the sample size
is $N=1377$.

We now discuss the results of applying Procedure \ref{proc:both-PC}.
Using dimension reduction in time only,  we obtained P-values
0.234 for $J=2$ (CPV = 92\%) and 0.220 for $J=3$ (CPV= 95\%).
Using the double dimension reduction, we obtained the following
values:
\begin{center}
\begin{tabular}{|c|c|c|}
\hline
&  \text{P--value} & CPV\\
\hline
$K=2, J=2$  & 0.272& 45\% \\
\hline
$K=3, J=3$  & 0.217 & 62\% \\
\hline
$K=4, J=4$ &  0.224& 67\% \\
\hline
$K=6, J=4$  & 0.223 & 80\% \\
\hline
$K=7, J=4$ &  0.221& 85\% \\
\hline
\end{tabular}
\end{center}
These remarkably similar P-values indicate that panels of
cumulative intraday return curves can in some cases be assumed to have
a separable covariance function. This could be useful for portfolio
managers as it indicates that they can exploit separability of the
data for more efficient modeling.

\medskip

\textcolor{black} {We conclude by noting that in practice it is important
to ensure that the time series forming the panel are at comparable scales.
This has been the case in our data examples, and will be the case if
the series are measurements of
 the same quantity and  are generated as a single group. If some of the series
are much more variable than the others, they may bias the test, and should
perhaps be considered separately.
}

\bigskip

\centerline{\small  ACKNOWLEDGMENTS}

\noindent This research has been
partially supported by the United States  National Science Foundation
grants at Colorado State University and Penn State University.

\bigskip

\centerline{\small SUPPORTING INFORMATION}

\noindent Additional Supporting Information may be found online in the
supporting information tab for this article.

\medskip

\renewcommand{\baselinestretch}{0.9}

\small

\bibliographystyle{oxford3}
\bibliography{sts}

\newpage

\renewcommand{\baselinestretch}{1.0}

\normalsize

\centerline{\Large \bf Supporting Information }

\appendix

\section{Proofs of the results of Section~\ref{s:da}} \label{s:p}

In \citetext{constantinou:2017}, the asymptotic
  distributions of the  test statistics were derived under the
  assumption of independent and identically distributed Gaussian data
  so that maximum likelihood estimators could be used to estimate the
  covariance and its separable analog. However, here we make no
  normality assumptions and we allow the sequence to be weakly
  dependent across $n$, thus entirely different proof techniques are
  employed.  In particular, we utilize multiple stochastic Taylor
  expansions to leverage the asymptotic normality of $\widehat
  C^{(h)}$ to derive the joint asymptotic distribution of $(\widehat
  C_1^{(h)}, \widehat C_2^{(h)}, \widehat C^{(h)})$ as well as the
  asymptotic behavior of our test statistics under both the null and
  alternative hypotheses.  These arguments become quite technical due
  to the fact that we are deriving asymptotic distributions of random
  operators.

\paragraph{Proof of Theorem \ref{thm:1}}
The starting point is the  asymptotic distribution of $\widehat{C}^{(h)}$.
It follows from Theorem 3 of \citetext{kokoszka:reimherr:2013}
that under  Assumption~\ref{a:main},
\begin{equation} \label{e:KR}
\sqrt{N}(\widehat{C}^{(h)} - C^{(h)})
\xrightarrow{\mathcal{L}} N ( 0, \boldsymbol{\Gamma}^{(h)}),
\end{equation}
 where $\boldsymbol{\Gamma}^{(h)}$ is given by
 \[
\boldsymbol{\Gamma}^{(h)} = \bR_0^{(h)} + \sum_{i=1}^\infty [\bR_i^{(h)} + (\bR_i^{(h)})^*]
\quad \text{with }
\bR_i^{(h)} =  \E [ ( \boldsymbol{X}_{1} \otimes \boldsymbol{X}_{1+h}- C^{(h)}) \otimes ( \boldsymbol{X}_{1+i} \otimes \boldsymbol{X}_{1+i+h} - {C}^{(h)})].
\]
Here $(\bR_i^{(h)})^*$ denotes the adjoint of $\bR_i^{(h)}$.  The
operator $\boldsymbol{\Gamma}^{(h)}$ is a Hilbert-Schmidt operator acting on
the space of Hilbert-Schmidt operators over $\mbH_1 \otimes \mbH_2$.
Since we have the asymptotic distribution of $\widehat{C}^{(h)}$, in
the following steps we use a one term Taylor expansion of the partial
trace operators to find the joint asymptotic distribution of
$\widehat{C}_1^{(h)},\widehat{C}_2^{(h)},\widehat{C}^{(h)}$.
Consider the operator:
\begin{equation*}
\bff (\widehat{C}^{(h)}) = \left( \begin{matrix}
\vspace{0.1 cm}
f_1(\widehat{C}^{(h)})\\
\vspace{0.1 cm}
f_2(\widehat{C}^{(h)})\\
\vspace{0.1 cm}
f_3(\widehat{C}^{(h)})
\end{matrix}
\right)
=\left( \begin{matrix}
\vspace{0.1 cm}
\frac {\Tr_2(\widehat{C}^{(h)})}{\Tr(\widehat{C}^{(h)})}\\
\vspace{0.1 cm}
 \Tr_1(\widehat{C}^{(h)})\\
 \vspace{0.1 cm}
 \widehat{C}^{(h)}
 \end{matrix}
\right)
=\left( \begin{matrix}
\vspace{0.1 cm}
\widehat{C}_1^{(h)}\\
\vspace{0.1 cm}
\widehat{C}_2^{(h)}\\
\vspace{0.1 cm}
\widehat{C}^{(h)}
 \end{matrix}
\right).
\end{equation*}
So $\bff(\hat C^{(h)})$ is an element of the Cartesian product space
$\mbH_1 \times \mbH_2 \times \mbH$ with $f_1:\mbH_1 \otimes \mbH_2
\rightarrow \mbH_1$, $f_2:\mbH_1\otimes \mbH_2 \rightarrow\mbH_2$,
$f_3:\mbH_1 \otimes
\mbH_2 \rightarrow \mbH_1\otimes \mbH_2$.
We will apply at Taylor expansion to $\bff(\hat C^{(h)})$ about the true parameter value $C^{(h)}$.  To do this, we require the Fr\'echet derivative of $\bff(\hat C^{(h)})$ which can be computed coordinate-wise as
\[
\frac{\partial}{\partial C^{(h)}}\bff =\left( \begin{matrix}
\frac {\partial{f_1}}{\partial{C^{(h)}}} \\
\frac {\partial{f_2}}{\partial{C^{(h)}}} \\
\frac {\partial{f_3}}{\partial{C^{(h)}}}
\end{matrix}\right).
\]
 Here $\partial f_i/\partial C^{(h)}$ denotes the Fr\'echet derivative of  $f_i$ with respect to  $C^{(h)}$.  Since $f_i$ is an operator, this means its derivative is a linear operator acting on the space of operators.
 Our goal is to use a Taylor expansion for Hilbert spaces to obtain the joint asymptotic distribution of
$\widehat{C}_1^{(h)},\widehat{C}_2^{(h)},\widehat{C}^{(h)}$.  We approximate $f_1(\widehat{C}^{(h)}), f_2(\widehat{C}^{(h)}), f_3(\widehat{C}^{(h)})$ by:
\begin{align*}
f_1(\widehat{C}^{(h)})& =  f_1(C^{(h)}) + \frac {\partial{f_1}}{\partial{C^{(h)}}}(\widehat{C}^{(h)} - C^{(h)}) + O_P(N^{-1}),\\
f_2(\widehat{C}^{(h)})& = f_2(C^{(h)})+ \frac {\partial{f_2}}{\partial{C^{(h)}}}(\widehat{C}^{(h)} - C^{(h)})+ O_P(N^{-1}),\\
f_3(\widehat{C}^{(h)})& = f_3(C^{(h)})+ \frac {\partial{f_3}}{\partial{C^{(h)}}}(\widehat{C}^{(h)} - C^{(h)})+ O_P(N^{-1}),
\end{align*}
where the last term is $O_P(N^{-1})$ because $\widehat{C}^{(h)} - C^{(h)} = O_P(N^{1/2})$.
In terms of the cartesian product form, this is equivalent  to:
\[
\bff (\widehat{C}^{(h)}) = \bff(C^{(h)}) +  \nabla\bff (C^{(h)}) (\widehat{C}^{(h)} - C^{(h)}) + O_P(N^{-1}).
\]
We therefore have that the variance operator of $\bff (\widehat{C}^{(h)})$ is asymptotically given by
\begin{align}\label{e:var}
\Var(\bff (\widehat{C}^{(h)})) &\approx  \Var \left(\bff(C^{(h)}) +  \nabla\bff (C^{(h)}) (\widehat{C}^{(h)} - C^{(h)}) \right)\\
&= \Var \left(\bff(C^{(h)}) +  \nabla\bff (C^{(h)})( \widehat{C}^{(h)}) - \nabla\bff (C^{(h)})(C^{(h)}) \right)\nonumber \\
&= \Var \left(\nabla\bff (C^{(h)}) (\widehat{C}^{(h)})\right)\nonumber \\
&= \nabla\bff (C^{(h)}) \left(\Var (\widehat{C}^{(h)}) (\nabla \bff (C^{(h)})^*) \right) \nonumber \\
&=\nabla\bff (C^{(h)}) \circ \boldsymbol{\Gamma}^{(h)} \circ \nabla \bff (C^{(h)})^* \nonumber \\
&:=\boldsymbol{W}^{(h)}\nonumber ,
\end{align}
where $(\cdot)^*$ denotes the adjoint operator.
We stress again that each term written above is a linear operator, and thus $\bW^{(h)}$ is actually a composition ($\circ$) of three linear operators.
This implies that the joint asymptotic distribution of $\widehat{C}_1^{(h)},\widehat{C}_2^{(h)},\widehat{C}^{(h)}$ is given by:
\begin{align*}
\sqrt{N}(\bff(\widehat{C}^{(h)}) - \bff(C^{(h)}))=\sqrt{N} \left( \begin{matrix}
\widehat{C}_1^{(h)}- C_1^{(h)} \\
\widehat{C}_2^{(h)}- C_2^{(h)} \\
\widehat{C}^{(h)} - C^{(h)}
\end{matrix}
\right)
\xrightarrow{\mathcal{L}} \mcN( {\bf 0}, \boldsymbol{W}^{(h)}).
\end{align*}

To complete the proof we need to find the
Fr\'echet derivatives. This turns out to be easier if we work with a
basis for the Hilbert spaces.  For example, the actions of a
continuous linear operator are completely determined by its actions on
individual basis elements.  Let $u_1,u_2, \dots$ be a basis for
$\mbH_1$ and $v_1,v_2,
\dots$ a basis for $\mbH_2$. Then a basis for $\mbH=\mbH_1\otimes
\mbH_2$ is given by $\{u_i \otimes v_j: i = 1,2,\dots , j =
1,2,\dots\}$.  Since $C_1^{(h)}: \mbH_1 \to \mbH_1$ is a compact operator, we can express it as
\[
C_1^{(h)} = \sum_{ik} C_{1;ik}^{(h)} u_i \otimes u_k,
\]
where $C_{1;ik}^{(h)} \in \mbR$ and $\sum_{ik} (C_{1;ik}^{{(h)}})^2 < \infty$.  Similarly we have that
\[
C_2^{(h)} = \sum_{jl} C_{2;jl}^{(h)} v_j \otimes v_l,
\]
and
\[
C^{(h)} = \sum_{ijkl} C_{ijkl}^{(h)} u_i \otimes v_j \otimes u_k \otimes v_l.
\]
These forms will be useful as we will be able to determine derivatives by taking derivatives with respect to the basis coordinate system.
In the following, $\delta_{ik}$ is the usual Kronecker delta.   We begin with $f_2(C^{(h)})$ as it is simpler than $f_1$.  Note that, by definition we have
\[
C_2^{(h)} = f_2(C^{(h)}) = \sum_{jl} \left( \sum_{i} C_{ijil}^{(h)} \right) v_{j} \otimes v_l,
\]
that is, we take the trace over the $u$ coordinates.  So we have hat $f_2$ is a linear mapping from $\mbH_1 \otimes \mbH_2 \otimes \mbH_1 \otimes \mbH_2 \to  \mbH_2 \otimes \mbH_2$.  If we take the derivative of this expression with respect to $C_{ijkl}^{(h)}$, then we get that
\[
\frac{\partial f_2(C^{(h)})}{\partial C_{ijkl}^{(h)}} = (v_j \otimes v_l) \delta_{ik}.
\]
If $i \ne k$, then $C_{ijkl}^{(h)}$ does not appear in the expression for $f_2(C^{(h)})$ and thus the derivative would be zero.  So we have that
\begin{align}
\frac{\partial f_2(C^{(h)})}{\partial C^{(h)}} = \sum_{ijkl} \frac{\partial f_2(C^{(h)})}{\partial C_{ijkl}^{(h)}} u_i \otimes v_j \otimes u_k \otimes v_l
= \sum_{ijkl} \delta_{ik} (v_j \otimes v_l) \otimes  (u_i \otimes v_j \otimes u_k \otimes v_l),
\label{e:Df2}
\end{align}
where again, this is interpreted as a linear operator from $\mbH_1 \otimes \mbH_2 \otimes \mbH_1 \otimes \mbH_2 \to  \mbH_2 \otimes \mbH_2$.  Note that the above operator is nearly the identity, e.g. $ \frac{\partial f_2(C^{(h)})}{\partial C^{(h)}} (x) = x$, but returns 0 for the off-diagonal $u$ coordinates.  We denote this operator as
\[
\bM_2^{(h)} = \sum_{ijkl} \delta_{ik} (v_j \otimes v_l) \otimes  (u_i \otimes v_j \otimes u_k \otimes v_l).
\]

The partial derivative of $f_1$ with respect to $C^{(h)}$ is a bit more complicated as it is a nonlinear function of $C^{(h)}$.  We can express $f_1$ as
\[
C_1^{(h)} = f_1(C^{(h)})
= \frac{\Tr_2(C^{(h)})}{\Tr(C^{(h)})}
=\sum_{ik} \frac{ \sum_{j} C_{ijkj}^{(h)}}{\sum_{i'j'} C_{i'j'i'j'}^{(h)}} u_{i} \otimes u_k.
\]
Again, taking the derivative with respect to the $C_{ijkl}^{(h)}$ coordinate, we get that
\[
\frac{\partial f_1(C^{(h)})}{\partial C_{ijkl}^{(h)}}
= \delta_{jl} \frac{\Tr(C^{(h)}) - \delta_{ik}\sum_{j'} C_{ij'kj'}^{(h)}}{\Tr(C^{(h)})^2} u_i \otimes u_k.
\]
Therefore we have that
\begin{align}
\frac{\partial f_1(C^{(h)})}{\partial C^{(h)}}
= \sum_{ijkl} \delta_{jl} \frac{\Tr(C^{(h)}) - \delta_{ik}\sum_{j'} C_{ij'kj'}^{(h)}}{\Tr(C^{(h)})^2} (u_i \otimes u_k) \otimes  (u_i \otimes v_j \otimes u_k \otimes v_l)
:= \bM^{(h)}_1.
\label{e:Df1}
\end{align}

Finally,  the partial derivative of $f_3$ with respect to $C^{(h)}$, i.e., the
partial derivative of $C$ with respect to $C$ is simply the identity, therefore
\begin{equation} \label{e:Df3}
\frac {\partial{f_3}}{\partial{C^{(h)}}}
=\sum_{ijkl} (u_i \otimes v_j \otimes u_k \otimes v_l) \otimes (u_i \otimes v_j \otimes u_k \otimes v_l) := \bM^{(h)}_3.
\end{equation}
By \eqref{e:KR}, \eqref{e:var},
 \eqref{e:Df1},  \eqref{e:Df2},
 and \eqref{e:Df3},
we obtain
\begin{align*}
\sqrt{N} \left( \begin{matrix}
\widehat{C}_1^{(h)}- C_1^{(h)} \\
\widehat{C}_2^{(h)}- C_2^{(h)} \\
\widehat{C} ^{(h)}- C^{(h)}
\end{matrix}
\right)
\xrightarrow{\mathcal{L}} N( {\bf 0}, \boldsymbol{W}^{(h)}),
\end{align*}
where $\boldsymbol{W}^{(h)}$ is given by
\begin{align}\label{e:wmatrix}
\bW ^{(h)} =
\left(
\begin{matrix}
\bM^{(h)}_1 \\
\bM^{(h)}_2 \\
\bM^{(h)}_3
\end{matrix}
\right)
\circ \bGa^{(h)}
\circ
\left(
\begin{matrix}
\bM^{(h)}_1 \\
\bM^{(h)}_2 \\
\bM^{(h)}_3
\end{matrix}
\right)^*.
\end{align}

\begin{remark} \label{r:W-blocks} The operator $\boldsymbol{W}^{(h)}$
has the following block structure form:
\begin{equation*}
 \boldsymbol{W} ^{(h)}=
\left (\begin{matrix}
 \boldsymbol{W}_{11}^{(h)} & \boldsymbol{W}_{12}^{(h)} &  \boldsymbol{W}_{13}^{(h)}\\
 \boldsymbol{W}_{21}^{(h)} & \boldsymbol{W}_{22}^{(h)} &  \boldsymbol{W}_{23}^{(h)}\\
 \boldsymbol{W}_{31}^{(h)} & \boldsymbol{W}_{32}^{(h)} &  \boldsymbol{W}_{33}^{(h)}
 \end{matrix} \right),
\end{equation*}
where
\begin{align*}
 \bW_{11}^{(h)} &= \bM_1^{(h)} \boldsymbol{\Gamma}^{(h)}(\bM^{(h)}_1)^*,\quad
 \bW_{12}^{(h)} =  \bM_1^{(h)}\boldsymbol{\Gamma}^{(h)} (\bM^{(h)}_2)^*, \quad
  \bW_{13}^{(h)} =\bM_1^{(h)}\boldsymbol{\Gamma}^{(h)} (\bM^{(h)}_3)^*\\
  \bW_{21}^{(h)} &=\bM_2^{(h)}  \boldsymbol{\Gamma}^{(h)} (\bM^{(h)}_1)^*, \quad
\bW_{22}^{(h)} =\bM_2^{(h)}   \boldsymbol{\Gamma}^{(h)}  (\bM^{(h)}_2)^*, \quad
 \bW_{23}^{(h)} = \bM_2^{(h)}  \boldsymbol{\Gamma}^{(h)} (\bM^{(h)}_3)^*\\
 \bW_{31}^{(h)} &=\bM_3^{(h)} \boldsymbol{\Gamma}^{(h)} (\bM^{(h)}_1)^*,\quad
  \bW_{32}^{(h)} =\bM_3^{(h)}  \boldsymbol{\Gamma}^{(h)} (\bM^{(h)}_2)^*, \quad
 \bW_{33}^{(h)} =\bM_3^{(h)}\boldsymbol{\Gamma}^{(h)}(\bM^{(h)}_3)^*=\boldsymbol{\Gamma}^{(h)}.
 \end{align*}
The operator  $\boldsymbol{W}_{11}^{(h)}$ is the covariance operator of
$\widehat{C}_1^{(h)}$ and $\boldsymbol{W}_{11}^{(h)} \in \cS(\cS(\mbH_1))$,
$\boldsymbol{W}_{12}^{(h)}$ is the covariance between $\widehat{C}_1^{(h)}$ and
$\widehat{C}_2^{(h)}$ and $\boldsymbol{W}_{12}^{(h)} \in \cS(\mbH_1) \otimes
\cS(\mbH_2)$, $\boldsymbol{W}_{13}^{(h)}$ is the covariance between
$\widehat{C}_1^{(h)}$ and $\widehat{C}^{(h)}$ and $\boldsymbol{W}_{13}^{(h)} \in
\cS(\mbH_1) \otimes \cS(\mbH_1 \otimes \mbH_2)$, $\boldsymbol{W}_{21}^{(h)}$
is the covariance between $\widehat{C}_2^{(h)}$ and $\widehat{C}_1^{(h)}$ and
$\boldsymbol{W}_{21}^{(h)} \in \cS(\mbH_2) \otimes \cS(\mbH_1)$,
$\boldsymbol{W}_{22}^{(h)}$ is the covariance operator of $\widehat{C}_2^{(h)}$
and $\boldsymbol{W}_{22}^{(h)} \in \cS(\cS(\mbH_2))$, $\boldsymbol{W}_{23}^{(h)}$
is the covariance between $\widehat{C}_2^{(h)}$ and $\widehat{C}^{(h)}$ and
$\boldsymbol{W}_{23}^{(h)}\in \cS(\mbH_2) \otimes \cS(\mbH_1 \otimes
\mbH_2)$, $\boldsymbol{W}_{31}^{(h)}$ is the covariance between
$\widehat{C}^{(h)}$ and $\widehat{C}_1^{(h)}$ and $\boldsymbol{W}_{31}^{(h)} \in
\cS(\mbH_1 \otimes \mbH_2) \otimes \cS(\mbH_1)$, $\boldsymbol{W}_{32}^{(h)}$
is the covariance between $\widehat{C}^{(h)}$ and $\widehat{C}_2^{(h)}$ and
$\boldsymbol{W}_{32}^{(h)}\in \cS(\mbH_1 \otimes \mbH_2) \otimes
\cS(\mbH_2)$ and finally $\boldsymbol{W}_{33}^{(h)}$ is the covariance
operator of $\widehat{C}^{(h)}$ and $\boldsymbol{W}_{33}^{(h)} \in \cS(\cS(\mbH_1
\otimes \mbH_2))$.
\end{remark}

\paragraph{Proof of Theorem \ref{thm:2}}
Since we have the joint asymptotic distribution of
$\widehat{C}_1^{(h)},\widehat{C}_2^{(h)},\widehat{C}^{(h)}$,  we can use the
delta method again to find the asymptotic distribution of
$\widehat{C}_1^{(h)}\widetilde\otimes \widehat{C}_2^{(h)} - \widehat{C}^{(h)}$ and in
particular, we can find the form of $\boldsymbol{Q}^{(h)}$, the asymptotic
covariance of $\widehat{C}_1^{(h)}\widetilde\otimes \widehat{C}_2^{(h)} -
\widehat{C}^{(h)}$.
Consider the function $g (\underaccent{\tilde}{C}^{(h)})= f_1(C^{(h)})
\widetilde \otimes f_2(C^{(h)}) - f_3(C^{(h)}) =
C_1^{(h)}\widetilde\otimes C_2^{(h)} -C^{(h)}$.  Using a one term
Taylor expansion we have that
\[
g(\widehat{\underaccent{\tilde}{C}}^{(h)}) \approx g (\underaccent{\tilde}{C}^{(h)}) + \frac {\partial(g (\underaccent{\tilde}{C}^{(h)}))}{\partial{C_1}^{(h)}}(\widehat{C}_1^{(h)} - C_1^{(h)}) + \frac {\partial(g (\underaccent{\tilde}{C}^{(h)}))}{\partial{C_2}^{(h)}}(\widehat{C}_2^{(h)} - C_2^{(h)})
+ \frac {\partial(g (\underaccent{\tilde}{C}^{(h)}))}{\partial{C}^{(h)}}(\widehat{C}^{(h)} - C^{(h)}),
\]
or equivalently
\[
g(\underaccent{\tilde}{\widehat{C}^{(h)}}) =  g (\underaccent{\tilde}{C}^{(h)}) + \nabla g (\underaccent{\tilde}{C}^{(h)})^* \left( \begin{matrix}
\widehat{C}_1^{(h)}- C_1^{(h)} \\
\widehat{C}_2^{(h)}- C_2^{(h)} \\
\widehat{C}^{(h)} - C^{(h)}
\end{matrix}
\right) + O_P(N^{-1}),
\quad \nabla g (\underaccent{\tilde}{C}^{(h)})^* = \left (\begin{matrix}
\frac {\partial(g (\underaccent{\tilde}{C}^{(h)}))}{\partial{C_1^{(h)}}},
\frac {\partial(g (\underaccent{\tilde}{C}^{(h)}))}{\partial{C_2^{(h)}}},
\frac {\partial(g (\underaccent{\tilde}{C}^{(h)}))}{\partial{C^{(h)}}}
\end{matrix} \right) ,
\]
which implies that the variance of $\widehat{C}_1^{(h)}\widetilde\otimes \widehat{C}_2^{(h)} - \widehat{C}^{(h)}$ is approximately:
\begin{align*}
\Var(\widehat{C}_1^{(h)}\widetilde\otimes \widehat{C}_2^{(h)} - \widehat{C}^{(h)})
\approx
&=\nabla g (\underaccent{\tilde}{C}^{(h)})^* \circ \boldsymbol{W}^{(h)} \circ \nabla g (\underaccent{\tilde}{C}^{(h)}) :=\boldsymbol{Q}^{(h)},
\end{align*}
and therefore the delta method implies that the asymptotic distribution of $\widehat{C}_1^{(h)}\widetilde\otimes \widehat{C}_2^{(h)} - \widehat{C}^{(h)}$ is given by:
\begin{align*}
\sqrt{N}(g(\widehat{\underaccent{\tilde}{C}}^{(h)}) - g (\underaccent{\tilde}{C}^{(h)}))=\sqrt{N} ( (\widehat{C}_1^{(h)}\widetilde \otimes \widehat{C}_2^{(h)}  - \widehat{C}^{(h)}) - (C_1^{(h)}\widetilde \otimes C_2^{(h)}  - C^{(h)}))\xrightarrow{\mathcal{L}} N ( 0, \boldsymbol{Q}^{(h)}).
\end{align*}

To complete the proof we need to find the partial derivatives. Taking
the derivative with respect to $C_{1}^{(h)}$ yields
\[
\frac {\partial(g (\underaccent{\tilde}{C}^{(h)}))}{\partial{C_1}^{(h)}} =\frac {\partial(C_1^{(h)}\widetilde\otimes C_2^{(h)} - C^{(h)})}{\partial{C_1}^{(h)}} = \boldsymbol{\mathcal{I}}_{4}\widetilde\otimes C_2^{(h)} = \boldsymbol{G}_1^{(h)},
\]
with respect to $C_2^{(h)}$
\[
\frac {\partial(g (\underaccent{\tilde}{C}^{(h)}))}{\partial{C_2}^{(h)}} =\frac {\partial(C_1^{(h)}\widetilde\otimes C_2^{(h)} - C^{(h)})}{\partial{C_2}^{(h)}} =C_1^{(h)}\widetilde\otimes  \boldsymbol{\mathcal{I}}_{4} =\boldsymbol{G}_2^{(h)},
\]
and with respect to $C^{(h)}$
\[
\frac {\partial(g (\underaccent{\tilde}{C}^{(h)}))}{\partial{C}^{(h)}} = \frac {\partial(C_1^{(h)}\widetilde\otimes C_2^{(h)} - C^{(h)})}{\partial{C}^{(h)}} = -\boldsymbol{\mathcal{I}}_8
\]
where $\boldsymbol{\mathcal{I}}_{4}$ and $\boldsymbol{\mathcal{I}}_8$ are the fourth and eighth order identity tensors.
Therefore by using the above partial derivatives we obtain the desired asymptotic distribution which is:
\begin{equation*}
\sqrt{N} ( (\widehat{C}_1^{(h)}\widetilde \otimes \widehat{C}_2^{(h)}  - \widehat{C}^{(h)}) - (C_1^{(h)}\widetilde \otimes C_2^{(h)}  - C^{(h)}))\xrightarrow{\mathcal{L}} N ( 0, \boldsymbol{Q}^{(h)}),
\end{equation*}
where $\boldsymbol{Q}^{(h)}$ is given by:
\begin{equation} \label{e:Qmatrix}
 \boldsymbol{Q}^{(h)} =
\left (\begin{matrix}
\boldsymbol{G}_1^{(h)}\\
\boldsymbol{G}_2^{(h)}\\
 -\boldsymbol{\mathcal{I}}_8
\end{matrix} \right)^*
\boldsymbol{W}^{(h)}
\left (\begin{matrix}
\boldsymbol{G}_1^{(h)}\\
\boldsymbol{G}_2^{(h)}\\
-\boldsymbol{\mathcal{I}}_8
 \end{matrix} \right).
\end{equation}

\begin{remark} \label{r:in:practice}
Here we provide more details on how $\bW^{(h)}_{KJ}$ and $\bQ^{(h)}_{KJ}$ are actually computed in practice.  Recall that
$\hat C^{(h)}_{KJ}$ is a $K \times J \times K \times J $ array, $\hat C_{1,K}^{(h)}$ is a $K \times K$ matrix, and $\hat C_{2,J}^{(h)}$ is a $J \times J$ matrix.  Also $\hat \bGa^{(h)}_{KJ}$ is a $K \times J \times K \times J \times K \times J \times K \times J $ array, and can be computed as described in Procedure~\ref{proc:both-PC}. To find $\bW^{(h)}_{KJ}$ we need to find first the array analogs of $\bM^{(h)}_1$, $\bM^{(h)}_2$, $\bM^{(h)}_3$, which we denote by $\bM^{(h)}_{1,KJ}$, $\bM^{(h)}_{2,KJ}$, $\bM^{(h)}_{3,KJ}$. $\bM^{(h)}_{2,KJ}$ is a $J \times J \times K \times J \times K \times J$ array with entries $M^{(h)}_{2,KJ}[j,l,i,j',k,l']=\delta_{jj'}\delta_{ll'}\delta_{ik}$, that is, it's an array with entries zeros and ones. $\bM^{(h)}_{3,KJ}$ is a $K \times J \times K \times J \times K \times J \times K \times J $ array with entries $M^{(h)}_{3,KJ}[i,j,k,l,i',j',k',l']=\delta_{ii'}\delta_{jj'}\delta_{kk'}\delta_{ll'}$, which is the eighth order identity array. Finally, $\bM^{(h)}_{1,KJ}$ is a $K \times K \times K \times J \times K \times J$ array. To find the entries of $\bM^{(h)}_{1,KJ}$ we first have to calculate $\boldsymbol{\Delta}$ and $\boldsymbol{\mathcal{I}}_{4}$. $\boldsymbol{\Delta}$ is a $K \times K \times K \times J \times K \times J$ array with entries $\Delta[i,k,i',j,k',l]=\delta_{ii'}\delta_{kk'}\delta_{jl}$, that is, it's an array with entries zeros and ones. $\boldsymbol{\mathcal{I}}_{4}$ is a $K \times J \times K \times J$ identity array, that is an array with entries $\mathcal{I}_{4}[i,j,k,l]=\delta_{ik}\delta_{jl}$. Then, we need to find $\Tr_2(\hat{C}^{(h)}_{KJ}) $, which is a $K\times K$ matrix with entries $\Tr_2(\hat{C}^{(h)}_{KJ})[i,k]= \sum_{j'} \hat{C}^{(h)}_{KJ,ij'kj'}$ and the scalar $\Tr (\hat{C}^{(h)}_{KJ})=\sum_{i'j'} \hat{C}^{(h)}_{KJ,i'j'i'j'}$. Combining $\boldsymbol{\Delta}$, $\boldsymbol{\mathcal{I}}_{4}$, $\Tr_2(\hat{C}^{(h)}_{KJ})$ and $\Tr (\hat{C}^{(h)}_{KJ})$ we can compute $\bM^{(h)}_{1,KJ}$ by:
\[
\bM^{(h)}_{1,KJ} = \frac{\boldsymbol{\Delta} \Tr (\hat{C}^{(h)}_{KJ}) - \Tr_2(\hat{C}^{(h)}_{KJ}) \otimes \boldsymbol{\mathcal{I}}_{4}}{(\Tr (\hat{C}^{(h)}_{KJ}))^2},
\]
where the tensor product can be easily implemented by using the {\tt R} package "tensorA" by
\citetext{boogaart:gerald:2007}.

Since we have $\bM^{(h)}_{1,KJ}$, $\bM^{(h)}_{2,KJ}$, $\bM^{(h)}_{3,KJ}$ we can compute $\bW^{(h)}_{KJ}$. Note that $\bW^{(h)}_{KJ}$ has a block structure of the following form:
\begin{equation}  \label{e:wmatrix:practice}
 \bW^{(h)}_{KJ} =
\left (\begin{matrix}
\bW^{(h)}_{11,K} & \bW^{(h)}_{12,KJ} &  \bW^{(h)}_{13,KJ}\\
 \bW^{(h)}_{21,KJ} & \bW^{(h)}_{22,J} &  \bW^{(h)}_{23,KJ}\\
 \bW^{(h)}_{31,KJ} & \bW^{(h)}_{32,KJ} &  \bW^{(h)}_{33,KJ}
 \end{matrix} \right),
\end{equation}
where
\begin{align*}
 \bW^{(h)}_{11,K} &= \bM^{(h)}_{1,KJ} \boldsymbol{\Gamma}^{(h)}_{KJ}(\bM^{(h)}_{1,KJ})^*,\quad
 \bW^{(h)}_{12,KJ} =  \bM^{(h)}_{1,KJ}\boldsymbol{\Gamma}^{(h)}_{KJ} (\bM^{(h)}_{2,KJ})^*, \quad
  \bW^{(h)}_{13,KJ} =\bM^{(h)}_{1,KJ}\boldsymbol{\Gamma}^{(h)}_{KJ} (\bM^{(h)}_{3,KJ})^*\\
  \bW^{(h)}_{21,KJ} &=\bM^{(h)}_{2,KJ}  \boldsymbol{\Gamma}^{(h)}_{KJ} (\bM^{(h)}_{1,KJ})^*, \quad
\bW^{(h)}_{22,J}=\bM^{(h)}_{2,KJ}  \boldsymbol{\Gamma}^{(h)}_{KJ} (\bM^{(h)}_{2,KJ})^*, \quad
 \bW^{(h)}_{23,KJ} = \bM^{(h)}_{2,KJ}  \boldsymbol{\Gamma}^{(h)}_{KJ} (\bM^{(h)}_{3,KJ})^*\\
 \bW^{(h)}_{31,KJ} &=\bM^{(h)}_{3,KJ} \boldsymbol{\Gamma}^{(h)}_{KJ} (\bM^{(h)}_{1,KJ})^*,\quad
  \bW^{(h)}_{32,KJ} =\bM^{(h)}_{3,KJ} \boldsymbol{\Gamma}^{(h)}_{KJ} (\bM^{(h)}_{2,KJ})^*, \quad
 \bW^{(h)}_{33,KJ} =\boldsymbol{\Gamma}^{(h)}_{KJ}.
\end{align*}
% \bW^{(h)}_{33,KJ} =\bM^{(h)}_{3,KJ}\boldsymbol{\Gamma}^{(h)}_{KJ}(\bM^{(h)}_{3,KJ})^*=\boldsymbol{\Gamma}^{(h)}_{KJ}.

$\bW^{(h)}_{11,K}$ is the variance--covariance array of
$\widehat{C}^{(h)}_{1,K}$ with dimensions $K \times K \times K \times K $,
$\bW^{(h)}_{12,KJ}$ is the covariance between $\widehat{C}^{(h)}_{1,K}$ and
$\widehat{C}^{(h)}_{2,J}$ with dimensions $K \times K \times J \times J $, $\bW^{(h)}_{13,KJ}$ is the covariance between
$\widehat{C}^{(h)}_{1,K}$ and $\widehat{C}^{(h)}_{KJ}$ with dimensions $K \times K \times K \times J \times K \times J $, $\bW^{(h)}_{21,KJ}$ is the covariance between $\widehat{C}^{(h)}_{2,J}$ and $\widehat{C}^{(h)}_{1,K}$ with dimensions
$J \times J \times K \times K$, $\bW^{(h)}_{22,J}$ is the variance--covariance array of $\widehat{C}^{(h)}_{2,J}$
with dimensions $J \times J \times J \times J$, $\bW^{(h)}_{23,KJ}$
is the covariance between $\widehat{C}^{(h)}_{2,J}$ and $\widehat{C}^{(h)}_{KJ}$ with dimensions $J \times J \times K \times J \times K \times J $, $\bW^{(h)}_{31,KJ}$ is the covariance between
$\widehat{C}^{(h)_{KJ}}$ and $\widehat{C}^{(h)}_{1,K}$ with dimensions $K \times J \times K \times J \times K \times K $,
$\bW^{(h)}_{32,KJ}$ is the covariance between $\widehat{C}^{(h)}_{KJ}$ and $\widehat{C}^{(h)}_{2,J}$ with dimensions $K \times J \times K \times J \times J \times J $ and finally $\bW^{(h)}_{33,KJ}$ is the variance--covariance
of $\widehat{C}^{(h)}_{KJ}$ which is $\boldsymbol{\Gamma}^{(h)}_{KJ}$.

To compute $\bQ^{(h)}_{KJ}$ we need to find the array analogs of the derivatives $\boldsymbol{G}^{(h)}_1$, $\boldsymbol{G}^{(h)}_2$ and $-\boldsymbol{\mathcal{I}}_8$, which we denote by $\boldsymbol{G}^{(h)}_{1,KJ}$, $\boldsymbol{G}^{(h)}_{2,KJ}$ and $-\boldsymbol{\mathcal{I}}_{8,KJ}$. First, notice that $-\boldsymbol{\mathcal{I}}_{8,KJ} = - \bM^{(h)}_{3,KJ}$. $\boldsymbol{G}^{(h)}_{1,KJ}$ is $K \times K \times K \times K \times J \times J $ array, which can be computed by the tensor product between the identity array of dimensions $K \times K \times K \times K$ and $\widehat{C}^{(h)}_{2,J}$. Similarly, $\boldsymbol{G}^{(h)}_{2,KJ}$ is $K \times K \times J \times J \times J \times J $ array, which can be computed by the tensor product between $\widehat{C}^{(h)}_{1,K}$ and the identity array of dimensions $J \times J \times J \times J$.

Since we have $\boldsymbol{G}^{(h)}_{1,KJ}$, $\boldsymbol{G}^{(h)}_{2,KJ}$ and $-\boldsymbol{\mathcal{I}}_{8,KJ}$ we can compute  $\bQ^{(h)}_{KJ}$. Note that $\bQ^{(h)}_{KJ}$ has the following form:
\begin{equation}  \label{e:qmatrix}
\bQ^{(h)}_{KJ} =\bQ^{(h)}_{1,KJ} +\bQ^{(h)}_{2,KJ} +  \bQ^{(h)}_{3,KJ}+
 \bQ^{(h)}_{4,KJ} +\bQ^{(h)}_{5,KJ} +  \bQ^{(h)}_{6,KJ}+
\bQ^{(h)}_{7,KJ} + \bQ^{(h)}_{8,KJ} +  \bQ^{(h)}_{9,KJ}
\end{equation}
where
\begin{align*}
 \boldsymbol{Q}^{(h)}_{1,KJ} &= (\boldsymbol{G}_{1,KJ}^{(h)})^* \boldsymbol{W}^{(h)}_{11,K} \boldsymbol{G}^{(h)}_{1,KJ}, \quad
 \boldsymbol{Q}^{(h)}_{2,KJ} =  (\boldsymbol{G}_{2,KJ}^{(h)})^* \boldsymbol{W}^{(h)}_{21,KJ} \boldsymbol{G}^{(h)}_{1,KJ}, \quad
 \boldsymbol{Q}^{(h)}_{3,KJ} =  -\boldsymbol{\mathcal{I}}_{8,KJ} \boldsymbol{W}^{(h)}_{31,KJ} \boldsymbol{G}^{(h)}_{1,KJ}; \\
 \boldsymbol{Q}^{(h)}_{4,KJ} &= (\boldsymbol{G}_{1,KJ}^{(h)})^* \boldsymbol{W}^{(h)}_{12,KJ} \boldsymbol{G}^{(h)}_{2,KJ}, \quad
 \boldsymbol{Q}^{(h)}_{5,KJ} = (\boldsymbol{G}_{2,KJ}^{(h)})^* \boldsymbol{W}^{(h)}_{22,K} \boldsymbol{G}^{(h)}_{2,KJ}, \quad
 \boldsymbol{Q}^{(h)}_{6,KJ} = -\boldsymbol{\mathcal{I}}_{8,KJ} \boldsymbol{W}^{(h)}_{32,KJ} \boldsymbol{G}^{(h)}_{2,KJ}; \\
 \boldsymbol{Q}^{(h)}_{7,KJ} &= - (\boldsymbol{G}_{1,KJ}^{(h)})^* \boldsymbol{W}^{(h)}_{13,KJ} \boldsymbol{\mathcal{I}}_{8,KJ}, \quad
 \boldsymbol{Q}^{(h)}_{8,KJ} = - (\boldsymbol{G}_{2,KJ}^{(h)})^* \boldsymbol{W}^{(h)}_{23,KJ} \boldsymbol{\mathcal{I}}_{8,KJ}, \quad
 \boldsymbol{Q}^{(h)}_{9,KJ} = \boldsymbol{W}^{(h)}_{33,KJ},
 %\boldsymbol{Q}^{(h)}_{9,KJ} = \boldsymbol{\mathcal{I}}_{8,KJ} \boldsymbol{W}^{(h)}_{33,KJ} \boldsymbol{\mathcal{I}}_{8,KJ} = \boldsymbol{W}^{(h)}_{33,KJ},
\end{align*}
where $\bQ^{(h)}_{KJ}$ and $\bQ^{(h)}_{i,KJ}$, i=1,\dots,9, are $K \times J \times K \times J \times K \times J \times K \times J $ arrays.
\end{remark}

\bigskip

\paragraph{Proof of Theorem \ref{t:under:alt}}
\begin{proof}

Let $T^*=\widehat{C}^{(h)}_1\widetilde \otimes \widehat{C}^{(h)}_2  - \widehat{C}^{(h)}$. Then we can write

\begin{align*}
\widehat{T}=N\|\widehat{C}^{(h)}_1\widetilde \otimes \widehat{C}^{(h)}_2  - \widehat{C}^{(h)}\|^2 &= N \langle \widehat{C}^{(h)}_1\widetilde \otimes \widehat{C}^{(h)}_2  - \widehat{C}^{(h)}, \widehat{C}^{(h)}_1\widetilde \otimes \widehat{C}^{(h)}_2  - \widehat{C}^{(h)} \rangle \\
&=N \langle T^{*} - \Delta + \Delta, T^{*} - \Delta + \Delta \rangle \\
&=N[\|\Delta\|^2 + 2 \langle T^{*} - \Delta , \Delta \rangle + \langle T^{*} - \Delta, T^{*} - \Delta \rangle]\\
&=N\|\Delta\|^2 + 2N^{1/2}\langle N^{1/2}(T^{*} - \Delta) ,
\Delta \rangle \\
& \ \ + \langle N^{1/2}(T^{*} - \Delta), N^{1/2}(T^{*} - \Delta) \rangle\\
&=N\|\Delta\|^2 +  O_P(N^{1/2}) +  O_P(1),
\end{align*}
and the claim follows.
\end{proof}

\section{Additional simulations} \label{s:sim:add}
In addition to the simulation results presented in  Section~\ref{s:sim},  we
consider here different values of the parameter $c$, i.e. the
parameter in the covariance function that controls separability. For
this scenario, we use $S=10$, $K=J=3$, $N=100,150,200$ and
$c=0,0.25,0.5,0.75,1$. The results are given in Table~\ref{tb:c-dep}.

\begin{table*}[ht]
\begin{center}
\begin{tabular}{|c|c|c|c|}
\hline
& $N=100$ & $N=150$ & $N=200$\\
\hline
$c=0$  & 4.8 & 5.5 & 6.1 \\
& 88\% & 89\% & 89\%\\
\hline
$c=0.25$  & 6.2 & 8.1 & 10.8\\
& 89\% & 90\% & 91\%\\
\hline
$c=0.5$ &  14.9 & 26.4 &35.6 \\
& 90\% & 91\% & 90\%\\
\hline
$c=0.75$  & 41.6 & 67.8 & 87.4 \\
& 90\% & 91\% & 92\%\\
\hline
$c=1$ &  73.3 & 95.2 & 99.5 \\
& 90\% & 91\% & 91\%\\
\hline
\end{tabular}
\end{center}
\vspace{0.3cm}
\caption{Dependence of the rejection rates
on the spatio--temporal interaction
parameter $c$ in \eqref{e:cov1}.  \label{tb:c-dep}}
\end{table*}

To supplement the results for the covariance function
\eqref{e:cov1}, we   consider here the following covariance function
\begin{equation} \label{e:cov2}
\sigma_{ss^\prime}(t,\tp)
=\frac{\sigma^2}{(a(t-t^\prime)^2+1)^{1/2}}
\exp\left(-\frac{b^2[|s - s^\prime|/(S-1)]^{2}}{(a(t - \tp)^2+1)^c}\right),
\end{equation}
which is a smoother version of the covariance function \eqref{e:cov1};
$|t-t^\prime|$ is replaced by $(t-t^\prime)^2$.  As a demonstration,
we set $a=3$, $b=2$, $\sigma^2=1$. As in Section~\ref{s:sim},
 we simulate the functions at
$T=50$ time points,  equally spaced on $[0,1]$, and $S=10$
coordinates in the panel. We consider only the case with dimension
reduction in both time and coordinates, under the null and alternative
hypothesis. The results are shown in Tables \ref{tb:null:nc}
and \ref{tb:al:nc}.

\begin{table*}[ht]\centering
{\footnotesize
\begin{tabular}{@{}c|cccccccccccc@{}}\toprule
 \multicolumn{1}{c}{}&\multicolumn{1}{c}{}&\multicolumn{3}{c}{$N=100$} & & \multicolumn{3}{c}{$N=150$} & & \multicolumn{3}{c}{$N=200$} \\
\cmidrule{3-5} \cmidrule{7-9} \cmidrule{11-13}
\multicolumn{1}{c}{}&&$J=2$ & $J=3$ & $J=4$ && $J=2$ & $J=3$ & $J=4$ && $J=2$ & $J=3$ & $J=4$ \\ \midrule

&$K=2$ & 6.5 &6.2 & 5.8 && 5.7& 7.5 & 5.5 && 6.0 & 7.1 & 6.8 \\
              & & (91\%)&(93\%)&(95\%)&&(93\%)&(94\%)&(95\%)&&(93\%)&(94\%)&(95\%)\\
$S=10$&$K=3$ & 6.6 & 4.7 & 7.0 && 7.4 & 6.2 & 5.5 && 7.3 & 5.9 & 5.8 \\
              & & (97\%)&(99\%)&(99\%)&&(97\%)&(99\%)&(99\%)&&(97\%)&(99\%)&(99\%)\\
&$K=4$ & 6.2 &6.0 & 5.8 && 7.6 & 6.9 & 6.6 && 5.4 & 5.8 & 5.4 \\
              && (96\%)&(99\%)&(99\%)&&(97\%)&(99\%)&(99\%)&&(97\%)&(99\%)&(99\%)\\
\bottomrule
\end{tabular}
\vspace{0.3cm}
\caption{Rejection rates under $H_0$ ($c=0$) for the covariance
function \eqref{e:cov2}. \label{tb:null:nc} }
}%footnotesize
\end{table*}

\begin{table*}[ht]\centering
{\footnotesize
\begin{tabular}{@{}c|cccccccccccc@{}}\toprule
 \multicolumn{1}{c}{}&\multicolumn{1}{c}{}&\multicolumn{3}{c}{$N=100$} & & \multicolumn{3}{c}{$N=150$} & & \multicolumn{3}{c}{$N=200$} \\
\cmidrule{3-5} \cmidrule{7-9} \cmidrule{11-13}
\multicolumn{1}{c}{}&&$J=2$ & $J=3$ & $J=4$ && $J=2$ & $J=3$ & $J=4$ && $J=2$ & $J=3$ & $J=4$ \\ \midrule

&$K=2$ & 29.8 &33.9 & 17.7 && 38.7 & 54.5 & 30.1 && 62.4 & 73.2 & 35.1 \\
              & & (88\%)&(85\%)&(83\%)&&(90\%)&(85\%)&(84\%)&&(87\%)&(85\%)&(78\%)\\
$S=10$&$K=3$ & 39.0 & 63.2 & 64.8 && 65.0 & 90.5 & 91.2 && 84.3 & 98.3 & 99.4 \\
              & & (97\%)&(98\%)&(97\%)&&(97\%)&(98\%)&(97\%)&&(96\%)&(98\%)&(96\%)\\
&$K=4$ & 37.8 & 63.8 & 65.4 && 67.7 & 89.7 & 93.6 && 87.5 & 98.2 & 99.8 \\
              && (97\%)&(99\%)&(99\%)&&(98\%)&(99\%)&(99\%)&&(98\%)&(99\%)&(99\%)\\
\bottomrule
\end{tabular}
\vspace{0.2cm}
\caption{Empirical power ($c=1$) for the covariance
function \eqref{e:cov2}. \label{tb:al:nc}}
}%footnotesize
\end{table*}

Finally, we check the performance of our test when $h=1$. For this
case, we simulate functional panels as the moving average process
\[
X_{ns}(t)=  e_{ns}(t)+e_{n-1s}(t),
\]
which is a 1-dependent functional time series. We generate $e_{ns}(t)$
as Gaussian processes with the following covariance function:
\begin{equation}\label{cov:h1}
\sigma_{ss^\prime}(t,\tp)
=\sigma^2 \exp\{-a[(t-\tp)^2+2\beta(t-\tp)(s-s^\prime)+(s-s^\prime)^2]\}.
\end{equation}
Clearly $\beta$ is the separability parameter, which takes values in
$[0,1)$. When $\beta=0$, we have a separable covariance. We set $a=3$
for our simulations.

For comparison, we add the simulations for $h=0$ by using the
covariance function \ref{cov:h1}. For $h=1$ the test tends to be
conservative, while for $h=0$ it overrejects. Consequently, the power
is higher for $h=0$.

\begin{table*}[ht]\centering
\begin{tabular}{@{}cccccccccccc@{}}\toprule
&\multicolumn{3}{c}{$N=100$} & & \multicolumn{3}{c}{$N=150$} & & \multicolumn{3}{c}{$N=200$} \\
\cmidrule{2-4} \cmidrule{6-8} \cmidrule{10-12}
&$J=2$ & $J=3$ & $J=4$ && $J=2$ & $J=3$ & $J=4$ && $J=2$ & $J=3$ & $J=4$ \\ \midrule

$S=4$ & 3.5 &1.8 & 1.5 && 4.0& 2.5 & 2.5 && 5.4 & 4.4& 4.2 \\
            & (94\%)&(99\%)&(100\%)&&(94\%)&(99\%)&(100\%)&&(95\%)&(99\%)&(100\%)\\
\bottomrule
\end{tabular}
\vspace{0.2cm}
\caption{Rejection rates under $H_0$ ($\beta=0$) for $h=1$
and  the covariances \eqref{cov:h1}.  \label{t:n:h=1}}
\end{table*}

\begin{table*}[ht]\centering
\begin{tabular}{@{}cccccccccccc@{}}\toprule
&\multicolumn{3}{c}{$N=100$} & & \multicolumn{3}{c}{$N=150$} & & \multicolumn{3}{c}{$N=200$} \\
\cmidrule{2-4} \cmidrule{6-8} \cmidrule{10-12}
&$J=2$ & $J=3$ & $J=4$ && $J=2$ & $J=3$ & $J=4$ && $J=2$ & $J=3$ & $J=4$ \\ \midrule

$S=4$ & 85.6 &89.5 & 89.5 && 99.4 & 99.8 & 99.8 && 100 & 100 & 100 \\
           & (94\%)&(99\%)&(100\%)&&(94\%)&(99\%)&(100\%)&&(94\%)&(99\%)&(100\%)\\
\bottomrule
\end{tabular}
\vspace{0.2cm}
\caption{Empirical power ($\beta=0.9$) for $h=1$
and  the covariances \eqref{cov:h1}. \label{t:al:h=1}}
\end{table*}

\begin{table*}[ht]\centering
\begin{tabular}{@{}cccccccccccc@{}}\toprule
&\multicolumn{3}{c}{$N=100$} & & \multicolumn{3}{c}{$N=150$} & & \multicolumn{3}{c}{$N=200$} \\
\cmidrule{2-4} \cmidrule{6-8} \cmidrule{10-12}
&$J=2$ & $J=3$ & $J=4$ && $J=2$ & $J=3$ & $J=4$ && $J=2$ & $J=3$ & $J=4$ \\ \midrule

$S=4$ & 6.5 &6.2 & 6.2 && 5.4& 5.5 & 5.3 && 6.4 & 5.5& 5.7 \\
            & (94\%)&(99\%)&(100\%)&&(94\%)&(99\%)&(100\%)&&(95\%)&(99\%)&(100\%)\\
\bottomrule
\end{tabular}
\vspace{0.2cm}
\caption{Rejection rates under $H_0$ ($\beta=0$)  for $h=0$
and  the covariances \eqref{cov:h1}. \label{t:n:h=0}}
\end{table*}

\begin{table*}[ht]\centering
\begin{tabular}{@{}cccccccccccc@{}}\toprule
&\multicolumn{3}{c}{$N=100$} & & \multicolumn{3}{c}{$N=150$} & & \multicolumn{3}{c}{$N=200$} \\
\cmidrule{2-4} \cmidrule{6-8} \cmidrule{10-12}
&$J=2$ & $J=3$ & $J=4$ && $J=2$ & $J=3$ & $J=4$ && $J=2$ & $J=3$ & $J=4$ \\ \midrule

$S=4$ & 100 & 100 & 100 && 100 & 100 & 100 && 100 & 100 & 100 \\
           & (94\%)&(99\%)&(100\%)&&(94\%)&(99\%)&(100\%)&&(94\%)&(99\%)&(100\%)\\
\bottomrule
\end{tabular}
\vspace{0.2cm}
\caption{Empirical power ($\beta=0.9$) for $h=0$
and  the covariances \eqref{cov:h1}. \label{t:al:h=0}}
\end{table*}

\end{document}